\documentclass{lms}
\usepackage{amsmath}
\usepackage{amstext}
\usepackage{amsfonts}
\usepackage{latexsym}
\usepackage{amscd}
\usepackage{xypic}
\newcommand{\la}{\ensuremath{\longrightarrow}}

\newcommand{\pone}{\ensuremath{\mathbb{P}^{1}}}

\newcommand{\sheaf}{\ensuremath{\mathcal{O}}}

\newtheorem{theorem}{Theorem}[section]
\newtheorem{lema}[theorem]{Lemma}
\newtheorem{proposition}[theorem]{Proposition}
\newtheorem{corollary}[theorem]{Corollary}
\newnumbered{claim}[theorem]{Claim}
\newnumbered{construction}[theorem]{Construction of contractions}
\newnumbered{assertion}{Assertion}    
\newnumbered{conjecture}[theorem]{Conjecture}  
\newnumbered{definition}[theorem]{Definition}
\newnumbered{hypothesis}[theorem]{Hypothesis}
\newnumbered{remark}[theorem]{Remark}
\newnumbered{note}[theorem]{Note}
\newnumbered{observation}[theorem]{Observation}
\newnumbered{problem}[theorem]{Problem}
\newnumbered{question}[theorem]{Question}
\newnumbered{algorithm}[theorem]{Algorithm}
\newnumbered{example}[theorem]{Example}
\newunnumbered{notation}[theorem]{Notation}

\title[D-Minimal Models, divisorial contractions]{Families of D-minimal models 
and applications to 3-fold divisorial contractions}

\author{Nikolaos Tziolas}

\classno{14E30 (primary), 14E35 (secondary).}

\begin{document}
\maketitle

\begin{abstract}
Let $X/T$ be a one parameter family of canonical 3-folds and let $D$ be a Weil divisor on it flat over $T$. 
We study the problem of when the $D_t$-minimal models of $X_t$ form a family and we obtain conditions 
for this to happen. As an application of this we classify 
terminal divisorial contractions $E \subset Y \la \Gamma \subset X$ contracting an 
irreducible surface $E$ onto the smooth curve $\Gamma$, in the case when the general section of $X$ through $\Gamma$ 
is a $D_5$ DuVal singularity.

\end{abstract}

\setcounter{section}{-1}
\section{Introduction}
Let $X$ be a 3-fold with canonical singularities and let $D$ be a Weil divisor on it. It is known~\cite{Kaw88} that 
the sheaf of graded $\sheaf_X$ algebras $R(X,D)=\oplus_{m\geq 0}\sheaf_X(mD)$ is finitely generated. 
Moreover this is equivalent to the existence of a birational morphism $Z\stackrel{g}{\la}X$ which is an isomorphism 
in codimension 1 and such that $g_{\ast}^{-1}D$ is $\mathbb{Q}$-Cartier and $g$-ample. Such a morphism is unique 
and we call it the $D$-minimal model of $X$. The process of making a Weil divisor 
$\mathbb{Q}$-Cartier is very important in the three dimensional minimal model program and it is of interest to get higher 
dimensional analogues.

In the first part of this paper 
we investigate whether this can be done in families of canonical 3-folds. More precisely, we study the following problem. 
Let $X/T$ be a one parameter family of canonical 3-folds and let $D$ be a Weil divisor in $X$ that is flat over $T$. 
Let $Z_t\stackrel{g_t}{\la}X_t$ be the $D_t$-minimal model of $X_t$, which exists as mentioned earlier.
We want to know if they form a family, i.e., if there is a morphism $Z/T\stackrel{g}{\la} X/T$ 
whose fibers over $T$ are the $D_t$-minimal models of $X_t$. 

A special case of this problem is the following. Let 
$Y/T \stackrel{f}{\la} X/T$ be a projective morphism between families of canonical 3-folds and 
assume that $Y_t \stackrel{f_t}{\la} X_t$ is a flopping contraction for all $t$. Do the flops of $f_t$ form a family? 
More precisely, is there a projective morphism 
$Y^{\prime}/T \stackrel{f^{\prime}}{\la} X/T$ such that $Y^{\prime}_t \stackrel{f_t^{\prime}}{\la}X_t$ 
is the flop of $f_t$, for all $t$? 
This problem is a special case of the previous one since $Y^{\prime}$ is the $D$-minimal model of $X$, where $D=f_{\ast}(-A)$, 
and $A$ is an $f$-ample divisor.
 
This last problem is directly related to the problem of whether the minimal model program runs in 
families of canonical 3-folds. It does work in families of terminal 3-folds~\cite{Ko-Mo92} but for families of canonical 
3-folds it is still unknown if it works. 

Example 2.16 shows that in general, $D$-minimal models do not form families. Therefore it is important to find conditions 
that need to be satisfied for the $D$-minimal models to form families. Theorems~\ref{crepant} and~\ref{deform} show that 
the $D_t$ minimal models form a family if $e(X_t;D_t)$, the number of crepant divisors of $X_t$ with center in $\text{supp}D_t$, 
is independent of $t$. 
Corollary 2.5 gives conditions for the flops to form families. However these conditions are probably too strong and the case of 
families of flops deserves a much more careful study.

As an application of the one parameter theory that was developed in the first part of this paper we study and classify certain 
3-fold divisorial contractions. Divisorial contractions is one of the two fundamental birational maps that appear 
in the three dimensional minimal model program. The other ones are the flips. 
Flips are known to exist and they have been classified~\cite{Ko-Mo92}. The structure of divisorial 
contractions is still not well understood and a detailed knowledge of them will have many applications and will 
contribute to the understanding of the birational structure of Mori fiber spaces through the Sarkisov program.

Let $E\subset Y \stackrel{f}{\la} \Gamma \subset X$ be a 3-fold divisorial contraction such that $X$ and $Y$ have terminal 
singularities. 
%$E$ is an irreducible divisor, $\dim \Gamma \leq 1$, $Y-E\cong X-\Gamma$, $-K_Y$ is $f$-ample and $X$ has terminal singularities. 
%It is called terminal if $Y$ has terminal singularities too. 
In the case when $\Gamma$ is a point, there are classification results due to Corti~\cite{Cor-Rei00} 
and Kawakita~\cite{Ka02},~\cite{Ka01}. 
The case when $\Gamma$ is a curve has been studied by Kawamata~\cite{Kaw94} and myself~\cite{Tzi02}. 
Suppose that $\Gamma$ is a curve. If a terminal 
divisorial contraction $f$ exists, then the general hyperplane section $S$ of $X$ through $\Gamma$ is DuVal~\cite{Ko-Mo92}. 
In~\cite{Tzi02} it is shown that given a smooth curve $\Gamma$ in a terminal 3-fold $X$, there is always a divisorial 
contraction $f$ with $Y$ canonical. In this paper we classify divisorial contractions 
in the case that the general hyperplane section $S$ of $X$ through $\Gamma$ is a $D_5$ DuVal singularity and $\Gamma$ is a 
smooth curve. I believe that the general $D_{2n+1}$ case can be treated in principle with same method that treats the $D_5$ case 
but it is much harder and the calculations involved are more complicated. 
The case when the general hyperplane section $S$ of $X$ through $\Gamma$ is a $D_{2n}$ 
DuVal singularity was studied in~\cite{Tzi02} and all such contractions were classified. Moreover, it was shown that the $D_{2n}$ 
cases are very different from the $D_{2n+1}$ cases.

Let $0\in X$ be an index one 3-fold terminal singularity and $\Gamma$ a smooth curve in $X$ 
passing through the singular point. A divisorial contraction $E\subset Y \stackrel{f}{\la} \Gamma \subset X$ 
contracting an irreducible divisor onto $\Gamma$ can be constructed as follows~\cite{Tzi02}. 
Let $g\colon W \la X$ be the blow up of $X$ along $\Gamma$. Then there are two $g$-exceptional divisors $E$ and $F\cong \mathbb{P}^2$. 
$E$ is a ruled surface over $\Gamma$ and $F$ is over the singularity of $X$. 
Let $h\colon Z \la W$ be the $E$-minimal model of $W$ which exists~\cite{Kaw88}. Then after running a 
suitable minimal model program we arrive to a $Z^{\prime}$ and the birational transform of $F$ in $Z^{\prime}$ can be contracted 
to get the required contraction $Y\stackrel{f}{\la} X$. 

In order to understand when a terminal contraction exists, it is important to understand the $E$-minimal model of $W$, $Z$. 
In general this is very difficult. In order to do this we obtain normal forms for the equations of $X$ and $\Gamma$ and we 
consider two cases. A general and a special one. In the general case we deform $\Gamma \subset X$ to a $\Gamma_0 \subset X_0$ 
having a simpler equation. Then we use the one parameter theory developed in the first part to show that in certain cases of interest, 
the divisorial contractions $E\subset Y\stackrel{f}{\la}\Gamma\subset X$ and 
$E_0 \subset Y_0 \stackrel{f_0}{\la}\Gamma_0 \subset X_0$ form a family. 
Then we explicitely construct the contraction $f_0$ and using results 
on deformations of terminal singularities we obtain information about $f$. In the special case we work directly with the equation of 
$X$. A complete classification of such contractions is given in Theorem~\ref{D5}.

This work has been done during my stay at the Max Planck Institute f\"{u}r Mathematik. 

Finally I would like to thank the referree for 
his many usefull comments and for helping me make this paper readable.

\section{Notation and terminology}

Let $X$ be a normal algebraic variety defined over the field of complex numbers $\mathbb{C}$, and 
$D$ a Weil divisor on it.

$D$ is called $\mathbb{Q}$-Cartier if there is an integer $m>0$ such that $mD$ is Cartier. 
We say that $X$ is $\mathbb{Q}$-factorial iff every Weil divisor on $X$ is $\mathbb{Q}$-Cartier.

Suppose that $K_X+\Delta$ is $\mathbb{Q}$-Cartier. 
Let $f \colon Y \la X$ is a birational morphism from a normal variety $Y$ and let $E_i \subset Y$ be the exceptional divisors. 
Then it is possible to write \[
K_Y+f_{\ast}^{-1}\Delta \equiv  f^{\ast}(K_X+\Delta) +\sum_{i}a(E_i,X,\Delta)E_i.\]
We call $a(E,X,\Delta)$ the discrepancy of $E$ with respect to $(X,\Delta)$. If $\Delta =0$, then we denote $a(E,X,0)$ by $a(E,X)$.
The pair $(X,\Delta)$ is called canonical (resp. terminal) 
iff $a(E,X,\Delta) \geq 0$ (resp. $a(E,X,\Delta)>0$ for all exceptional divisors $E$~\cite{Rei87}.

An $f$-exceptional divisor $E \subset Y$ 
is called $f$-crepant if $a(E,X,\Delta)=0$. The birational morphism $f$ is called crepant iff all $f$-exceptional divisors 
are $f$-crepant. Let $Z\subset X$ be a closed subscheme. Then we define the number of crepant divisors 
with center in $Z$ by \[
e(X;Z)=\#\{E;a(E,X)=0 \;\; \text{and} \;\; \text{center}_X E \subset Z\}.\]
By a slight abuse of notation we set $e(X;D)=e(X;\text{supp}(D))$, for a Weil divisor $D$.
If $Z=0$ then we denote $e(X;0)$ by $e(X)$.
If  $X$ is canonical then it is known that $e(X;Z)$ is finite~\cite{Ko-Mo98}. 

\begin{definition} 
The $D$-minimal model of $X$ is a projective morphism $f\colon Y\la X$ such that $Y$ is normal, the exceptional set of $g$ 
has codimension at least 2, $D^{\prime}=f_{\ast}^{-1}D$ is $\mathbb{Q}$-Cartier and $f$-ample over $X$. If it exists then it 
is unique.  
\end{definition}
It is known that the existence of the $D$-minimal model is equivalent to the finite generation of the sheaf of algebras 
$R(X,D)=\oplus_{m\geq 0}\sheaf_X(mD)$~\cite{Kaw88}. 

Many important constructions in birational geometry are related to $D$-minimal models. The flop is one of the most important ones. 
A flopping contraction is a projective birational morphism $f\colon Y \la X$ to a normal variety $X$ 
such that the exceptional set of $f$ has codimension at least 2 in $Y$ and $K_Y$ is numerically $f$-trivial.
If $D_Y$ is a $\mathbb{Q}$-Cartier divisor on $Y$ such that $-(K_Y+D_Y)$ is $f$-ample, 
then a projective  morphism $f^{\prime}\colon Y^{\prime} \la X$ such that the exceptional set of $f^{\prime}$ has codimension 
at least 2 and that $K_{Y^{\prime}}+D_{Y^{\prime}}$ is $f^{\prime}$-ample is called the $D$-flop of $f$. 
In other words, it is the $D$-minimal 
model of $X$, where $D=f_{\ast}(D_Y)$. Its existence is fundamental in the context of the minimal model program. 

If $X$ is a 3-fold with canonical singularities, then $R(X,D)$ is 
finitely generated and the $D$-minimal model exists. The higher dimensional case is one of the most important 
open problems of higher dimensional birational geometry. 

A $\mathbb{Q}$-factorialization of $X$ is a projective birational morphism $f\colon Y \la X$ such that $Y$ is normal 
and $\mathbb{Q}$-factorial and $f$ is an isomorphism in codimension 1. If $X$ is a canonical 3-fold then a 
$\mathbb{Q}$-factorialization as above exists~\cite{Kaw88}.

Let $X$ be a canonical n-fold. A terminalization of $X$ is a crepant projective birational morphism $f\colon Y \la X$ 
such that $Y$ has terminal only singularities. If $n=3$ then a terminalization exists~\cite{Rei83}.

Let $f \colon Y \la X$ be a birational morphism between normal varieties and let $D$ be a Weil divisor in $X$. Then we denote 
the birational transform $f_{\ast}^{-1}D$ of $D$ in $Y$ by $D^Y$.

Let $\mathcal{X}\la T$ be a flat family of algebraic varieties. A family of Weil divisors over $T$ is a Weil divisor 
$D=\sum_ia_iD_i$ on $\mathcal{X}$ such that $D_i$ is flat over $T$ for all $T$. 

\section{Families of $D$-minimal models}
In order to decide whether the minimal model program runs in families of canonical 3-folds, it is important 
to investigate if $D$-minimal models can be constructed in families. More precisely we want to study the following problem.
 
\noindent\textbf{Problem:} Let $X \la T$ be a family of canonical 3-folds and let $D$ be a family of divisors in $X$. 
It is known that the $D_t$-minimal model of $X_t$ exists~\cite{Kaw88} for all $t \in T$. 
Do they form an algebraic family? i.e., 
Does there exists a morpism $f \colon Y \la X$ such that $f_t \colon Y_t \la X_t$ is the $D_t$ minimal model of $X_t$?

In particular, we want to know if canonical flops exist in families. More precisely.

\noindent\textbf{Problem:} Let $ Y/T \stackrel{f}{\la} X/T$ 
be a projective morphism between families of canonical 3-folds such that $Y_0 \stackrel{f_0}{\la} X_0$ is a flopping contraction. 
Does there exist a projective morphism $Y^{\prime}\stackrel{f^{\prime}}{\la} X$ such that 
$Y^{\prime}_0 \stackrel{f_0^{\prime}}{\la} X_0$ is the flop of $f_0$? 

Let $A$ be an $f$-ample divisor and $D=f_{\ast}(-A)$. Then $Y^{\prime}$ is the $D$-minimal model of $X$ and therefore 
this problem is a special case of the problem concering families of $D$-minimal models.  

From now on and unless otherwise said, $T$ will denote a smooth curve. Then the total spaces $X$ and $Y$ are 4-folds. 
Recently Shokurov showed that canonical 4-fold log flips exist~\cite{Sho02} and therefore the minimal model program works 
in dimension four. Let $f \colon Y/T \la X/T$ be as above.
Then the 4-fold flop of $f$ exists. Suppose it is $f^{\prime} \colon Y^{\prime} \la X$. 
However this does not necessarily mean that $Y^{\prime}_0 \la X_0$ is the flop of $f_0$. The reason is that even though 
$f$ is a small contraction, it is possible that $f_t$ is small for all $t \neq 0$, but $f_0$ is divisorial. 
The next example shows that this is indeed possible.

\begin{example}
Let $X_0 \subset \mathbb{C}^4$ be the 3-fold defined by the equation $xy+yu+xu=0$, where $x,\; y,\; z, \; u$ 
are the coordinates in $\mathbb{C}^4$. It's singularities 
are cDV points and its singular locus is the line $L_0: x=y=u=0$. Therefore $X_0$ is canonical. 
Let $X \la \mathbb{C}^1$ be the one parameter smoothing of $X_0$ given by $xy+yu+xu+ut=0$, where $t$ is the parameter. 
Let $Z$ in $X$ be the family of divisors given by $x=u=0$, and let $f\colon Y \la X$ be the blow up of $X$ along $Z$. 
Then $f_t \colon Y_t \la X_t$ is the blow up of $X_t$ along the divisor $Z_t$ given by $x=u=0$. 
If $t \neq 0$, then $X_t$ is smooth and 
therefore $f_t$ is an isomorphism. On the other hand, $X_0$ is singular along $L_0 \subset Z_0$ and hence $f_0$ is divisorial 
contracting an irreducible divisor onto $L_0$.
\end{example}  
The case of a family of terminal threefolds is particularly nice and it is possible to give a complete answer. 
\begin{theorem}\label{terminal}
Let $X \la T$ be a one dimensional family of terminal 3-folds. 
Let $D$ be a family of Weil divisors in $X$ over $T$. Then there exists a morphism $f \colon Y \la X$, 
such that $f_t\colon Y_t \la X_t$ is the $D_t$-minimal model of $X_t$ for all $t$.
\end{theorem}
\begin{proof}
First I claim that $X$ has terminal singularities. This follows from the following results.
\begin{theorem}[~\cite{Kaw99}]
Let $f\colon X\la T$ be a flat morphism from a germ of an algebraic vatiety to a germ of a smooth curve. Assume that the central fiber $X_0=f^{-1}(P)$ has only canonical singularities. 
Then so has the total space $X$ as well as any fiber $X_t$ of $f$. Moreover, the pair $(X,X_0)$ is canonical.  
\end{theorem}
It is also known~\cite{Nak98} that if $X_0$ is terminal, then so is $X$ and $X_t$. 
An immediate consequence of the above Theorem is that if $g \colon W \la X$ is a resolution of $X$, then there is no $g$-exceptional crepant divisor with center in $X_0$. Hence 
all crepant exceptional divisors of $X$ dominate $T$.

Now let $g \colon W \la X $
be a resolution of $X$ and let $E_i$ be the exceptional divisors. Run a $(W,0)$ MMP over $X$. Since $X$ is terminal, we arrive at a 
$Y^{\prime} \stackrel{g^{\prime}}{\la} X$ over $T$, such that $g^{\prime}$ is an isomorphism in codimension 1 and $Y^{\prime}$ is 
$\mathbb{Q}$-factorial. Therefore $D^{\prime}={g^{\prime}_{\ast}}^{-1}D$ is $\mathbb{Q}$-Cartier, and hence so is 
$D^{\prime}_t$. To make it $g^{\prime}_t$-ample, run a $(Y^{\prime}_t, -\epsilon D^{\prime}_t)$ MMP over $X_t$. 
By~\cite[Prop. 11.4, Theorem 11.10]{Ko-Mo92}, and since $X$ is terminal, this can be done in families. Therefore we obtain 
$Y \stackrel{g}{\la} X$ with the required properties.
\end{proof}

The case of families of canonical 3-folds is much more complicated. In general it is not possible to construct $D$-minimal models 
in families as shown by example~\ref{ex1}. 
The number of crepant divisors of the members of the family is important as the next Theorem shows.

\begin{theorem}\label{crepant}
Let $X \stackrel{\sigma}{\la} T$ be a proper family of canonical threefolds over a smooth curve $T$. Let $D$ in $X$ be a family of divisors in $X$ over $T$ 
and let $0\in T$ be a closed point. Then \[
e(X_0) \geq e(X_t) \]
for $t$ in a small neighborhood of $0 \in T$. Moreover,
\begin{enumerate}
\item If $e(X_0)=e(X_t)$ for all $t$, then after a finite base change there is a morphism $Y \stackrel{g}{\la} X$, with $Y$ $\mathbb{Q}$-factorial and terminal, 
such that $Y_t \stackrel{g_t}{\la} X_t$ is a terminalization of $X_t$, for all $t\in T$. 
\item If $e(X_0;D_0)=e(X_t;D_t)$ for all $t$, then there is a morphism $Y \stackrel{g}{\la} X$, such that 
$Y_t \stackrel{g_t}{\la} X_t$ is the $D_t$-minimal model of $X_t$.
\end{enumerate}
 \end{theorem}

We are now able to get some information about the existence of flops in families.
\begin{corollary}\label{flops}
Let $Y/T \stackrel{f}{\la} X/T$ be a morphism between one parameter families of canonical 3-folds. Suppose that 
$f_0 \colon Y_0 \la X_0$ is a flopping contraction. If $e(X_0)=e(X_t)$ for all $t$, then there is a morphism 
$Y^{\prime}/T \stackrel{f^{\prime}}{\la} X/T$ such that $f_0^{\prime}\colon Y_0^{\prime}\la X_0$ is the flop of $f_0$. 
\end{corollary}
\begin{proof}
Let $D_Y$ be a divisor in $Y$ such that $-D_Y$ is $f$-ample. Let $D=f_{\ast}(D_Y)$. 
Then $-D_{Y,t}$ is $f_t$-ample for all $t$ and therefore $Y_t$ is the $(-D_t)$-minimal model of $X_t$, and the flop 
of $f_t$ is the $D_t$-minimal model of $X_t$. Since $e(X_0)=e(X_t)$ for all $t$, then from Theorem~\ref{crepant} follows 
that there exist a morphism $f^{\prime}\colon Y^{\prime} \la X$ such that $f^{\prime}_t\colon Y^{\prime}_t \la X_t$ 
is the flop of $f_t$.
\end{proof}
\begin{remark} The condition that $e(X_t)$ is constant in the family is probably too strong and  
the $D$-minimal models may form a family even though this condition does not hold. 
For example, with notation as in the previous Corollary, suppose that $f_0$ 
 contracts a chain of rational curves to a canonical double point on $X_0$ but $X_t$ is only terminal for $t\neq 0$. 
Then the the flop of $f_0$ can be  constructed in families as 
Propostion~\ref{doublepoints} shows even though $e(X_0) >e(X_t)=0$, for $t\neq 0$. 
\end{remark}

\begin{proof}[of Theorem~\ref{crepant}]
The claimed inequality between the number of crepant divisors of the central and general fiber is independent of base change and so we will perform it when necessary. 
First I claim that up to a finite base change we can assume that $e(X_t)=e(X)$, for general $t\in T$. Indeed, let $g\colon W \la X$ be a log resolution of $X$ such that all the $g$-crepant 
divisors are smooth. Let $E$ be a $g$-crepant 
divisor. From the previous theorem it follows that the center of $E$ is not contained in $X_0$ and since there are finitely many crepant divisors, then after removing finitely many points from $T$ 
we can assume that all the $g$-crepant divisors dominate $T$. By generic smoothness, $E_t$ and $W_t$ are smooth, for general $t$,  
and hence $g_t \colon W_t\la X_t$ is a resolution of $X_t$. Then by adjunction it follows that $E_t$ is $g_t$-crepant. The problem here is that $E_t$ may have more than one connected components and 
so in general we only get that $e(X_t)\geq e(X)$. However, I claim that after a finite base change, $E_t$ is irreducible and smooth, for all irreducible crepant exceptional divisors $E$ of $X$, 
and therefore $e(X_t)=e(X)$. So, let $E$ be an irreducible $g$-crepant divisor of $X$ and let $E\stackrel{h}{\la} D \stackrel{\tau}{\la} T$ be the Stein factorization of 
$f=\tau \circ g \colon E \la T$, 
where $D=\mathrm{Spec}(f_{\ast}\sheaf_E)$. Then $h$ has connected fibers and $\tau$ is finite. We now make a base change with $D\la T$. So, let $X_D=X\times_T D$, 
$W_D = W\times_T D$ and $E_D=E\times_T D$. Note that by the previous theorem, $X_D$ is also canonical since all the fibers $X_{D,d}$ are canonical. 
Also $W_{D,d}$ is smooth for general $d \in D$. Moreover, by the universal property of fiber product we see that there is 
an embedding $E\subset E_d$, and by construction, $E\la D$ has connected fibers. Hence $E_d$ is irreducible, smooth and crepant for $X_{D,d}=X_t$, for general $t=\tau (d)$. 
Since $D$ may not be normal, make another base change with its normalization $\overline{D}$. Repeat this 
process for any crepant divisor $F$ of $W_{{\overline{D}}}$ such that $F_d$ is not irreducible. Since the crepant divisors are at most $e(X_t)$, this process ends with a family 
$X^{\prime} \la T^{\prime}$, with $T^{\prime}$ smooth, such that there is a log resolution $f^{\prime} \colon W^{\prime} \la X^{\prime}$, such that if $E$ is any $f^{\prime}$-exceptional 
crepant divisor, $E_t$ is smooth and irreducible and hence $e(X^{\prime}_t)=e(X^{\prime})$, for general $t$.  

We may also assume that $X$ is $\mathbb{Q}$-factorial. If this is not the case, then let $f \colon Y\la X$ be a $\mathbb{Q}$-factorialization, which exists by the MMP in dimension 4. 
Then $e(X)=e(Y)$. Moreover, 
$Y_t\la X_t$ is an isomorphism in codimension 1 for general $t$ and hence $e(Y_t)=e(X_t)=e(X)=e(Y)$. The central fiber contraction $Y_0 \la X_0$ may be divisorial, but in any case it is crepant 
and hence $Y_0$ is normal and canonical and $e(Y_0)\leq e(X_0)$. Hence in addition to our hypothesis we may also assume that $X$ is $\mathbb{Q}$-factorial and $e(X_t)=e(X)$. 

Let $n=e(X)$ and let $E_1,\ldots,E_n$ be the crepant exceptional divisors of $X$. Then by standard MMP arguments, we may extract them from $X$ with a series of crepant morphisms 
\begin{equation}
Y=X_n \stackrel{f_n}{\la} X_{n-1}\stackrel{f_{n-1}}{\la} X_{n-2} \la \cdots \la X_1 \stackrel{f_1}{\la} X
\end{equation}
where $f_i\colon X_{i}\la X_{i-1}$ is crepant and its exceptional set is $E_i$. Then by Theorem 2.3 the centers of these divisors are not contained in $X_0$.  
Hence $X_{i,0}$ is irreducible and $E_i \cap X_{i,0}$ is a divisor. Moreover, I claim than $X_{i,0}$ is normal and canonical. Indeed, inductively it easily follows 
that $K_{X_i}+X_{i,0}=f_i^{\ast}(K_{{X_{i-1}}}+X_{i-1,0})$ and hence since $(X,X_0)$ canonical, $(X_i,X_{i,0})$ is canonical as well and therefore $X_{i,0}$ is normal. Moreover, 
by adjunction it follows that $K_{{X_{i,0}}}=f_i^{\ast}K_{{X_{i-1,0}}}$ and hence $X_{i,0}$ is canonical and 
$E_i \cdot X_{i,0}$ is a crepant divisor for $X_{i-1,0}$ (which may be reducible). Therefore $e(X_0)\geq n=e(X_t)$. 

Suppose now that $e(X_0) = e(X_t)$. Let $f \colon Y\la X$ be the composition of the maps $f_i$ in $(2.1)$ above. Then by its construction, $Y$ is a $\mathbb{Q}$-factorial terminal 4-fold and since 
$e(Y_t)=0$ for general $t$, $Y_t$ is terminal as well for general $t$. Moreover, by the above discussion it follows that $Y_0$ is irreducible and $e(Y_0)=0$, and hence $Y_0$ 
is terminal too. Now $Y\la X$ satisfies all the conditions of Theorem~\ref{crepant}.1. 

Now suppose that $D \subset X$ is a family of divisors such that $e(X_0;D_0)=e(X_t;D_t)$, for all $t$. Let $Z \stackrel{g}{\la} X$ be the $D$-minimal model of $X$, which exists 
by the MMP in dimension 4. Then $Z_t \la X_t$ is an isomorphism in codimension 1 for general $t$ and therefore $e(Z_t)=e(X_t)$. I now claim that $Z_0 \la X_0$ is also an isomorphism 
in codimension 1 and hence it is also the $D_0$-minimal model of $X_0$. Suppose not. Let $D^{\prime}=f^{-1}_{\ast}D$. Then by the definition of $D$-minimal models, 
$-D^{\prime}$ is $g$-ample and hence if $Z_0 \la X_0$ is divisorial, then the center of any $g_0$-exceptional divisor is contained in $D_0$. Therefore, 
$e(D^{\prime}_0; Z_0)<e(D_0;X_0)=e(D_t;X_t)=e(D^{\prime}_t; Z_t)$, which is impossible from the first part of the proof.

\end{proof}
\begin{remark*} 
The condition $e(X_0)=e(X_t)$ is not sufficient for the existence of a morphism $g \colon Z\la X$ such that $g_t \colon Z_t \la X_t$ 
is a $\mathbb{Q}$-factorialization of $X_t$, for all $t$, as it was mistakenly claimed in~\cite{Tzi05}. The reason is that there may be divisors in $X_0$ that do not deform with $X_0$, 
i.e., they do not extend to a divisor in $X$. For example, let $X=(xy-zu+t=0) \subset \mathbb{C}^4$. Then $X_0=(t=0)$ is the ordinary double point $xy-zu=0$, and $X_t$ for $t\neq 0$ is smooth. 
$X_0$ is not $\mathbb{Q}$-factorial and there is no morphism $g\colon Y \la X$ such that $g_0 \colon Y_0\la X_0$ is a $\mathbb{Q}$-factorialization of $X_0$ because 
if there was such a morphism $g$, then $g$ would be an isomorphism in codimension 1 which is impossible since $X$ itself is smooth. 
\end{remark*}

Unfortunately, the proof of Theorem~\ref{crepant} is using the existence of the MMP in dimension 4. It would be desirable 
to have a proof of it that does not involve the MMP, like in dimension 3, and constructive if possible. 
In general I do not know how to do this but in the next special case that is of interest in the study of divisorial contractions, 
it can be done. Moreover, in this case there is no restriction on the dimension of the base $T$. 
 
\begin{theorem}\label{deform}
Let $X \la T$ be a family of index 1 canonical 3-folds, with $T$ smooth. Let $Z \subset X$ be a smooth effective divisor in $X$ 
flat over $T$. Moreover, assume that 
\begin{enumerate}
\item $X_t^{\text{sing}}\cap Z_t$ is a smooth curve $L_t$. Let $S_t=\text{Spec} \sheaf_{X_t,L_t}$. This is a DuVal surface 
singularity defined over $\mathbb{C}(t)$. 
$Z_t|_{S_t}$ is a line through the origin such that the extended dual graph of $S_t$ and $Z_t|_{S_t}$ is independent of $t$. 
This means that $S_t$ are isomorphic 
singularities and that the position of $Z_t|_{S_t}$ in the fundamental cycle does not change with $t$.
\item $X_0$ is either cDV along $Z_0$, or there is a smooth curve $L \subset X$ that dominates $T$ such that $L \cap X_t$ is a 
reduced point $P_t$ and 
$\mathrm{mult}_L X=2,3$ or $\geq 3$. This means that $P_t \in X_t$ is elliptic with constant $k=\mathrm{mult}_L X$. Moreover, 
the singularities of the weighted blow up  $B^w_{P_0} X_0$ over $P_0$ on the birational transform $\tilde{Z_t}$
are finitely many cDV points, 
and $w=(1,1,1,1)$, $(3,2,1,1,)$ or $(2,1,1,1)$ if $k=1$, $2$, or $\geq 3$ respectively.
\end{enumerate}
Then there is a birational morphism $Y \la X$ flat over $T$, such that $Y_t \la X_t$ is the $Z_t$-minimal model of $X_t$, 
for all $t$.
\end{theorem}
\begin{remark}
The complicated assumptions of the previous theorem easily imply  that $e(X_0;Z_0)=e(X_t;Z_t)$, for all $t$, and therefore 
if $\dim T =1$, Theorem~\ref{deform} 
follows from Theorem~\ref{crepant}. However, its proof does not use the MMP in dimension 4 and it is constructive. 
\end{remark}
The proof of Theorem~\ref{deform} is given at the end of section 2. 
The idea is the following. By performing a series of well chosen blow ups we arrive at a crepant morphism $h\colon W/T \la X/T$ 
suct that $W_t$ is irreducible and that $(h_t)_{\ast}^{-1}Z_t$ is $\mathbb{Q}$-Cartier. We then run a relative MMP in families 
on $W$ and arrive at a morphism $Y \la X$ with the properties claimed in the Theorem. 

The difficulty here is that in general $W_t$ is only canonical and therefore we want to run a MMP in families of canonical 3-folds.  
For families of terminal 3-folds it is known to work~\cite{Ko-Mo92}. The only operations of the MMP that is not known 
if they can be performed in families of canonical 3-folds is the flop and flip. Since we are always concerned with crepant morphisms, 
only the flop is needed in our case. 

In what follows we obtain some information about cases when flops can be performed in families. 
We also study the behaviour of the property of being $\mathbb{Q}$-factorial in families of canonical 3-folds, and we also 
obtain a few more technical results that are needed for the proof of Theorem~\ref{deform}.

In the terminal case, the property of being $\mathbb{Q}$-factorial is stable under deformations as the next Theorem shows.
\begin{theorem}[{\cite[Theorem 12.1.10]{Ko-Mo92}}]
Let $X \stackrel{f}{\la} T$ be a 1-dimensional family of terminal 3-folds, and let $Z \subset X$ be proper 
and flat over $T$.  
Let $D$ be a divisor on $X$ such that $D_0=D|_{X_0}$ 
is $\mathbb{Q}$-Cartier along $Z_0$. Then there is a neighborhood $V$ of $0 \in T$ such that $D$ is $\mathbb{Q}$-Cartier 
along $f^{-1}(V) \cap Z$.
\end{theorem}
In particular it follows that $D_t$ is $\mathbb{Q}$-Cartier along $Z_t$ for $t$ in a neighborhood of $0 \in T$. 

Unfortunately the previous result is no longer true for families of canonical 3-folds and one must be very careful 
when dealing with such families. The following example shows that it is possible that $D_0$ is $\mathbb{Q}$-Cartier, 
but $D_t$ is not for all $t \neq 0$.

\begin{example} We will now construct a one parameter family of canonical 3-folds $X \la T$, a closed subscheme $Z \subset X$ 
proper and flat over $T$ and a divisor $D \subset X$, such that $D_0$ is $\mathbb{Q}$-Cartier along $Z_0$, but $D_t$ is not 
$\mathbb{Q}$-Cartier along $Z_t$ for all $t \neq 0$.

Start with the family given by $Y: x^2+y^2z+xz^5+t u^3 +u^{10}=0$, with $t$ the parameter. Let $\Gamma \subset Y$ 
be given by $x=y=u=0$. Now let $W=B_{\Gamma}Y \stackrel{g}{\la}Y$ be the blow up of $Y$ along $\Gamma$. There are two exceptional 
divisors $E$ and $F$. $E$ is a \pone\ bundle over $\Gamma$, and $F$ a $\mathbb{P}^2$ bundle over the line $x=y=z=u=0$. 
Now let $X=B_E W \stackrel{h}{\la} W$ be the blow up of $W$ along $E$. Let $E^X$, $F^X$ be the birational transforms 
of $E$ and $F$ in $X$. Then, $F^X$ is proper over $T$, and I claim that $E^X_0$ is $\mathbb{Q}$-Cartier along $F^X_0$, but 
$E^X_t$ is not $\mathbb{Q}$-Cartier for all $t \neq 0$.  

For $t= 0$, $Y_0$ is given by $x^2+y^2z+xz^5+u^{10}=0$. Now explicit calculations of the blow up 
in exactly the same way as those in the proof 
of Theorem~\ref{D5}, show that $E_0+5 F_0$ is Cartier and $E^X_0$ is $\mathbb{Q}$-Cartier along $F^X_0$. 
Moreover, $E_0$ is not Cartier. Suppose it was. Then take a general fiber of $g_0$, $l \subset E_0$. This is rationally equivalent to 
a line $l^{\prime} \subset F_0 =\mathbb{P}^2$. Then $-1=E_0 \cdot l = E_0 \cdot l^{\prime} >0 $.

Let us check now what happens in nearby fibers.
For $t\neq 0$, $Y_t$ is given by $x^2+y^2z+xz^5+u^3=0$. Calculate the blow up of $\Gamma$ now. 
In the affine chart $u \neq 0$ we have the coordinates $x=xu$, $y=yu$. Then $W_0$ is given by $x^2u+y^2zu+xz^5+ u^2=0$, 
and $E_0, \; F_0$ are given by $x=u=0$ and $z=u=0$ respectively. Now blow up $E_0$ to get $X_0$. In the affine chart 
$u=ux$, $X_0$ is given by $x^2u+y^2zu+z^5+u^2x=0$, $E^{X_0}_0$ by $x=y^2u+z^4=0$, $F^{X_0}_0$ by $z=u=0$, and 
the $h$-exceptional divisor $B$ by $x=z=0$. Now it easy to see that $5F^{X_0}_0$ is Cartier and that $X_0$ is smooth away 
from $E^{X_0}_0$. But then, if $E^{X_0}_0$ was $\mathbb{Q}$-Cartier along $F^{X_0}_0$, then it would in fact be $\mathbb{Q}$-Cartier. 
But this would imply that $E_0$ is also $\mathbb{Q}$-Cartier as well which is not true. 
\end{example}
We want to run a MMP in families of canonical 3-folds. To do this it is important to study the problem of whether canonical 
flops form families. 
This is true for terminal flops~\cite{Ko-Mo92}. For canonical I do not know the answer in general but Theorem~\ref{crepant} 
suggests that in general it should be no. It is then very interesting to find conditions under which canonical flops 
can be constructed in families. Corollary~\ref{flops} gives some conditions in this direction but they are most likely too 
strong. If the singularities that we are studying are double points then flops can be constructed in families as the 
next propistion shows.

\begin{proposition}\label{doublepoints}
Let $f \colon Y/T \la X/T $ be a morphism between families of $n-$folds over $T$. Suppose that $f_0\colon Y_0 \la X_0$ is a
flopping contraction and that $X_0$ has only hypersurface double points. 
Then the flop $f_0^{\prime} \colon Y_0^{\prime} \la X_0$ of $f_0$ exists and there is a projective morphism 
$f^{\prime}: Y^{\prime} \la X$, such that $f_0^{\prime}: Y_0^{\prime} \la X_0$ is the flop of $f_0$. 
\end{proposition}
\begin{proof}
By using the Weierstrass preparation theorem we may write $0 \in X_0$ as \[
x^2 + f(y_1, \ldots , y_n)=0. \]
Therefore, $X$ can be written as \[
x^2+f(y_1,\ldots,y_n,t_1,\ldots,t_m)=0.\]
where $t_1,\ldots,t_m$ are local analytic coordinates of $T$. Now we proceed as in the 3-fold terminal case~\cite{Ko-Mo92}. 
If $-D$ is $f$-ample, then to prove the existence of the flop for $f$ it suffices 
to show that $R(X,D)=\oplus_d \sheaf_X (df_{\ast}D)$ is finitely generated. $X$ has an involution $i$. Then $D+i^{\ast}D=0$. Since 
$R(X,-D)$ is finitely generated, then so is $R(X,D)$. Therefore there is a 
morphism $f^{\prime}\colon Y^{\prime} \la X$ which is the flop of $f$. Moreover, $Y^{\prime}=Y$, $f^{\prime}=i \circ f$ and 
the involution is also in every fiber. Therefore the flop exists in families. 
\end{proof}
\begin{proposition}
Let $f \colon Y/T \la X/T$ be a morphism between two families of canonical 3-folds over $T$. Suppose that $f_0 \colon Y_0 \la X_0$ is 
a flopping contraction. Suppose that there is a morphism $g \colon Z \la Y$ over $T$, 
such that $Z_0 \la Y_0$ is a $\mathbb{Q}$-factorial 
terminalization of $X_0$. Then there is a morphism $f^{\prime} \colon Y^{\prime} \la X$, such that $f^{\prime}_0 \colon Y^{\prime}_0 \la X_0$, 
is the flop of $f_0$. 
\end{proposition}
Note that the conditions of this proposition are not always satisfied as shown by example~\ref{ex1}.
\begin{proof}
Let $f_0^{\prime} \colon Y_0^{\prime} \la X_0 $ be the flop of $f_0$. Let $g_0^{\prime} \colon Z_0^{\prime} \la Y_0^{\prime}$ 
be a $\mathbb{Q}$-factorial terminalization of $Y_0^{\prime}$. $Z_0$, $Z_0^{\prime}$ are birational and nef over $X_0$. Therefore, 
$\phi_0 \colon Z_0 \dasharrow Z_0^{\prime}$ is an isomorphism in codimension 1 and a composition of terminal flops~\cite{Kaw88}. 
Terminal flops exist in families~\cite{Ko-Mo92} and therefore $\phi_0$ extends to a map $\phi \colon Z \dasharrow Z^{\prime}$ over $T$. 
Then by~\cite{Ma-Ro71}, there is a morphism $g^{\prime} \colon Z^{\prime} \la Y^{\prime}$ over $X$, extending $g_0^{\prime}$. 
Hence $f^{\prime} \colon Y^{\prime} \la X$ is the required flop.
\end{proof}
The following lemmas are needed for the proof of Theorem~\ref{deform}.
\begin{lema}
Let $Z \subset X$ be a smooth divisor in a normal variety $X$. Suppose that $X$ has canonical hypersurface 
singularities only, and that every 
irreducible component of $X^{sing} \cap Z $ is smooth. Let $Y=B_Z X \stackrel{f}{\la} X$ be the blow up of $X$ along $Z$. 
Let $Z^{\prime}=f_{\ast}^{-1}Z$, be the birational transform of $Z$ in $Y$. Then
\begin{enumerate}
\item Y is normal and has hypersurface singularities only.
\item $K_Y = f^{\ast}K_X$.
\item Every irreducible component $E_i$ of the $f$-exceptional divisor $E$ is smooth.
\item Every component of $Y^{sing}\cap Z^{\prime}$ is smooth.
\end{enumerate}
\end{lema}
\begin{proof}
Let $n=\dim X$. The results are all local around the singularities of $X$ and therefore 
we may assume that $ X \subset \mathbb{C}^{n+1}$. Let $f \colon B_Z \mathbb{C}^{n+1} \la \mathbb{C}^{n+1}$ 
be the blow up of $\mathbb{C}^{n+1}$ along $Z$. Then $B_Z \mathbb{C}^{n+1}$ is smooth and $Y=B_Z X \subset B_Z \mathbb{C}^{n+1}$ 
has codimension 1 and hence $Y$ has hypersurface singularities only. Let $F$ be the $f$-exceptional divisor. Then by adjunction, 
$K_{B_Z \mathbb{C}^{n+1}}=f^{\ast}K_{\mathbb{C}^{n+1}} + F$, and $f^{\ast}X=Y+F$. Therefore \[
K_Y=(K_{B_Z \mathbb{C}^{n+1}}+Y)|_Y=f^{\ast}K_X. \]
Since $\mathbb{C}^{n+1}$ is smooth and $X$ canonical, the proof of ~\cite[Theorem 5.34]{Ko-Mo98} show that the pair $(\mathbb{C}^{n+1},X)$ 
is also canonical. By adjunction we see that \[
K_{B_Z \mathbb{C}^{n+1}}+Y=f^{\ast}(K_{\mathbb{C}^{n+1}} +X).\]
Therefore the pair $(\mathbb{C}^{n+1},Y)$ is also canonical. Hence by~\cite[Proposition 5.51]{Ko-Mo98}, $Y$ is normal. Finally, since 
every irreducible component of $X^{sing} \cap Z $ is smooth and $f^{-1}(z)=\mathbb{P}^1$, for all $z \in Z$, it follows that every 
irreducible component of $F\cap Y$ is also smooth.

Now let $L$ be a component of the singular locus of $X$ that lie on $Z$. Then this is either a point or a smooth curve. 
If it is a point then there is nothing to show. So assume it is a smooth curve. At the generic point of $L$, $X$ is a surface DuVal 
singularity, say $X_L$. Let $U \la X_L$ be the minimal resolution and $E_i$, $i=1 , \ldots ,n$ the exceptional curves. 
$Z$ corresponds to a line through the singularity. We will now consider cases with respect to the type of the singularity 
$X_L$ and the position of $Z$ in the dual graph. 

\noindent \textit{Case 1:} $X_L$ is an $A_n$ singularity and $Z$ intersects $E_k$ in the dual graph.

\noindent \textit{Subcase 1:} $(k,n)=(m+1,2m+1)$, for some $m$. 
Then by~\cite{Tzi02}, the $f$-exceptional divisor $F$ is $E_{m+1}$, and hence 
again by~\cite{Tzi02}, $Z^{\prime}$ is smooth over the generic point of $L$, and then $Z^{\prime}\cap Y^{sing}$ over $L$ is  
a finite set of points and perhaps a fiber of $f$ which is a \pone. 
Note that in this case it is possible that a component of $Y^{sing}$ is a 
singular curve, but it does not lie on $Z^{\prime}$. 

\noindent \textit{Subcase 2:} $(k,n)\neq (m+1,2m+1)$, for any $m$. In this case, $Y_L$ has two singular point, $P \in Z^{\prime}$ 
and $Q$. Therefore $Y^{sing}\cap Z^{\prime}$ consists of two curves over $L$ and possibly some fibers of $f$. It is now 
clear that all these curves must be smooth. 

\noindent \textit{Case 2:} $X_L$ is either $D_n$ or $E_i$, $i=6,\;7,\;8$ type DuVal singularity. 
These cases are treated exactly similar as 
the $A_n$ cases and we will not do them.
\end{proof}
\begin{lema}
Let $(0 \in X)$ be a canonical 3-fold singularity such that the general hyperplane section $H$ of $X$ through $0$ is an elliptic 
surface singularity with invariant $k=1$, $2$ or $3$. Let $0 \in Z \subset X$ be a smooth divisor. Let $Y \stackrel{f}{\la} X$ 
be either the blow up of $X$ at $0$ 
if $k=3$, or the $(3,2,1,1)$ or $(2,1,1,1)$ weighted blow up of $X$ at $0$ in the case that $k=2$ or $1$. 
Then the birational transform $Z^{\prime}$ of $Z$ in $Y$ is smooth and $Y^{sing}\cap Z $ is a smooth rational curve or a point or 
empty.
\end{lema}
\begin{proof}
\setcounter{equation}{0}
Suppose that $k=3$. Then $Y$ is just the blow up of $X$ at $0$ and therefore $Z^{\prime}=B_0Z$ is smooth.

Now suppose that $k=2$ or $1$. Then $\text{mult}_0X=2$~\cite[Theorem 4.57]{Ko-Mo98}. 
Let $Z$ be given by $x=y=0$. Then $X$ is given by an equation of the form \[
xf(x,y,z,t)+yg(x,y,z,t)=0.\]
Since $t=0$ is elliptic, we can write the equation of $X$ as \[
x^2+xf(x,y,z,t)+yg(x,y,z,t)=0.\]
By using the Weierstrass preparation theorem, the equation of $X$ can be written as 
\begin{equation}
x^2+2x\phi (y,z,t)+yh(y,z,t)=0.
\end{equation}
Moreover, $\phi (y,z,t), \; h(y,z,t) \in (y,z,t)^2$, because otherwise $0 \in X$ is $cA_n$. Eliminate $x$ from $(1)$. 
Then the equation of $X$ becomes 
\begin{equation}
F(x,y,z,t)=x^2+yh(y,z,t)-\phi^2(y,z,t)=0,
\end{equation}
and $Z$ is given by $x-\phi(y,z,t)=y=0$. Moreover from the previous discussion, $\phi(y,z,t)$, $h(y,z,t) \in (y,z,t)^2$. 
Now consider cases with respect to $h$.

\noindent\textit{Case 1:} $\text{mult}_0h(y,z,t)\geq 3$. Then assign weights to $x$, $y$, $z$ and $t$ as follows. 
Let $w(x)=2$ and $w(y)=w(z)=w(t)=1$. 
Let $Y=B_0^wX \stackrel{f}{\la} X$ be the $(2,1,1,1)$ weighted blow up of $X$. 
We want to understand the birational transform $Z^{\prime}$ 
of $Z$ in $Y$.

\noindent\textit{Claim:} $Z^{\prime}=B_0Z$ and therefore $Z$ is smooth. By the definition of the weighted blow up, \[
Z^{\prime}=\text{Proj} \oplus_{n \geq 0} \frac{\textstyle m^{w}(n)+I_Z}{\textstyle I_Z},\]
and $I_Z=(x-\phi(y,z,t),y)$. We may look at $Z$ as given by $x-\phi_1(z,t)=0$ in $\mathbb{C}^3_{x,z,t}$, 
where $\phi_1(z,t)=\phi(0,z,t)$. 
Note that $\phi_1(z,t) \neq 0$, since otherwise $\phi(y,z,t)=y \lambda (y,z,t)$ and hence $2Z$ given by $x^2=y=0$ will be Cartier. 
Then $\text{mult}_0\phi_1(z,t) \geq 2$. Now 
\begin{gather*}
\frac{\textstyle m^{w}(1)+I_Z}{\textstyle I_Z}=(\overline{x},\overline{z},\overline{t})=(\overline{z},\overline{t},
\phi_1(\overline{z},\overline{t}))= 
(\overline{z},\overline{t})=\frac{\textstyle m +I_Z}{\textstyle I_Z}
\end{gather*}
In general, 
\begin{gather*}
\frac{\textstyle m^{w}(n)+I_Z}{\textstyle I_Z}=(\overline{x}^i\overline{t}^j\overline{z}^k; 2i+j+k \geq n)=\\
(\phi_1(\overline{z},\overline{t})^i\overline{t}^j\overline{z}^k; 2i+j+k \geq n)=(\overline{t}^j\overline{z}^k; j+k\geq n)=
\frac{\textstyle m^n +I_Z}{\textstyle I_Z}
\end{gather*}
Therefore, $Z^{\prime}=B_0Z$, and hence $Z^{\prime}$ is smooth. Note that this is the case that $k=2$. 

\textit{Case 2:} $\text{mult}_0h(y,z,t)=2$, and $k=1$. 
Let $h_2(y,z,t)$ be the degree 2 part of $h(y,z,t)$. Since the section $t=0$ is elliptic, the cubic term of $F(x,y,z,0)$ must be a cube 
because otherwise the section $t=0$ is DuVal. Therefore, \[
h_2(y,z,t)=ay^2+tl(y,z,t),\]
where $l(y,z,t)$ is linear in $y$, $z$, $t$. Therefore the equation of $X$ becomes 
\begin{equation}
x^2+ay^3+ytl(y,z,t)+yh_{\geq 4}(y,z,t)-\phi^2(y,z,t)=0.
\end{equation}
If $l(y,z,t) \neq 0$, then one of $y^2z$, $yzt$, $yt^2$, appears in the above equation and the section $z=t$ is given by $x^2+\psi(y,t)=0$, 
and $\psi_3(y,t)$ is not a cube and hence it must be DuVal. Therefore $l(y,z,t)=0$. Therefore the equation of $X$ becomes 
\begin{equation}
x^2+y^3+yh_{\geq 4}(y,z,t)-\phi^2(y,z,t)=0,
\end{equation} 
where $\text{mult}_0 \phi \geq 2$. Moreover, $Z$ is given by $x-\phi(y,z,t)=y=0$, and as before, $\phi(0,z,t)\neq 0$. 
By using the Weierstrass 
preparation theorem it is possible to write 
\begin{gather*}
x^2+y^3+yh_{\geq 4}(y,z,t)-\phi^2(y,z,t)=\\
x^2+u_1[y^3+\alpha (z,t) y^2 + \beta (z,t) y + \gamma (z,t)]= \\
u_1[(\frac{x}{\sqrt{u_1}})^2 + y^3 + \alpha (z,t) y^2 +\beta (z,t) y + \gamma(z,t)]
\end{gather*}
From the above equations, it is clear that $\gamma(z,t)=unit \cdot \phi^2_{\geq 2}(0,z,t)=-\psi_{\geq 2}^2(z,t)$. 
Therefore $Z$ is given by $x-\psi_{\geq 2}(z,t)=y=0$. Now eliminate $y^2$. The equation of $X$ can be written \[
x^2+y^3+3\alpha(z,t)y^2+\beta(z,t)y-\psi^2_{\geq 2}(z,t)=0.\]
The change of variables $y \mapsto y - \alpha(z,t)$ makes the equation of $X$ 
\begin{equation}
x^2+y^3+(y-\alpha(z,t))\delta(z,t)-\alpha^3(z,t)-\psi^2(z,t)=0,
\end{equation}
where $\delta(z,t)=\beta(z,t)-3\alpha^2(z,t)$. $Z$ is given by $x-\psi_{\geq 2}(z,t)=y-\alpha (z,t)=0$. Now it is easy to see that 
$\text{mult}_0(\alpha^3(z,t)+\psi^2(z,t)) \geq 6$. Therefore, 
$\text{mult}_0(\psi) \geq 3$ and $\text{mult}_0(\alpha) \geq 2$. Now an argument as in case 1 shows that 
$Z^{\prime}=B_0Z$ and therefore it is smooth.

Finally it remains to justify the statement about the singular locus of $Y$. The singularities of $Y$ on $Z^{\prime}$ over $0$ are on 
the intersection of $Z^{\prime}$ and the exceptional divisors. But since $f$ is a certain weighted blow up, this intersection 
is a smooth rational curve. 
\end{proof}
The next elementary result is needed for the proof of Theorem~\ref{deform}.
\begin{lema}
Let $T$ be a smooth curve and $Y/T \stackrel{f}{\la} X/T$ be a morphism over $T$. Assume that $X_0$ is a canonical 3-fold. 
Suppose that $K_Y=f^{\ast}K_X$. Then $\dim f^{-1}(x) \leq 2$ for 
any $x \in X_0$.
\end{lema}
\begin{proof}
As before we see that the pair $(X,X_0)$ is canonical. Therefore, for any $x_0 \in X_0$, $x_0 \in X$ is a terminal singularity. 
If there is a divisor $E \subset f^{-1}(x_0)$, then $E$ is crepant over $x_0$ which is impossible. Therefore $\dim f^{-1}(x_0) \leq 2$.
\end{proof}

\subsection{Proof of Theorem~\ref{deform}}

Suppose that $X_0$ is given by $F(x,y,z,w)=0$ in $\mathbb{C}^4$. Then $X$ is given by $F(x,y,z,w)+t\Psi(x,y,z,w)=0$ in $\mathbb{C}^5$. 
Let $X_1\stackrel{f_1}{\la} X$ be the $(1,1,1,1,1)$, or $(2,1,1,1,1)$, or $(3,2,1,1,1)$ weighted blow up of $X$ along $L$. 
Consider the $(1,1,1,1,1)$ blow up case first. Let $\tilde{\mathbb{C}}=B_0\mathbb{C}^5$. Then 
\begin{gather*}
K_{\tilde{\mathbb{C}}}=f_1^{\ast}K_{\mathbb{C}^5}+3E \\
f_1^{\ast}X=Y+3E.
\end{gather*}
Therefore \[
K_{X_1}=f_1^{\ast}K_X .\]
By lemma 2.15, it follows that $\dim f_1^{-1}(x)\leq 2$ for any $x \in X_0$. Therefore $X_1 \la T$ is also a 
family of canonical 3-folds. 
Moreover, $e(X_1) < e(X)$. By lemma 2.14, $Z_1=(f_1)_{\ast}^{-1}Z$ is smooth and $X_1^{sing}\cap (Y_1)_0$ is a union of smooth rational 
curves. Let $X_2 \stackrel{f_2}{\la} X_1$ be the blow up of $X_1$ along $Z_1$. Then by lemma 2.13, $K_{X_2}=f_2^{\ast}K_{X_1}$, $X_2$ is 
normal and $e(X_2) < e(X_1)$. Moreover, the conditions of part 1 of the theorem guarantee that there is a smooth $f_2$-exceptional 
divisor (over the generic point of $L_0$ it is just  the blow up of a line through a DuVal singularity). 
Moreover, again by lemma 2.13, every component of the singular locus of $X_2$ that lies on $Z_2$ is smooth. 
If it is just an isolated set of points then stop. If not, then blow up the $f_2$-exceptional divisor that they lie on, 
which as mentioned is smooth, to get a crepant morphism $X_3 \stackrel{f_3}{\la} X_2$, 
with $e(X_3) <  e(X_2)$. Continue this process of blowing up exceptional divisors to get a 
sequence \[
\cdots \la X_k \stackrel{f_k}{\la} X_{k-1} \la \cdots \stackrel{f_3}{\la} X_2 \stackrel{f_2}{\la} X_1 \stackrel{f_1}{\la} X,\]
such that \[
e(X_1)>e(X_2)> \cdots e(X_k) > \cdots \]
Since there are finitely many crepant divisors, this process must stop and there is an $n>0$ such that either $e(X_n)=0$ 
and hence $X_n$ is terminal, or by lemma 2.13 there are finitely many only singularities on $Z_n$.
Moreover, $g_n \colon X_n \la X$ is crepant and there are finitely many singularities only on every $g_n$-exceptional curve. 
Let $E_i$ be the $g_n$-exceptional divisors which by lemma 2.13 are all smooth. 
As we mentioned above, $(X_n)_0$ has only isolated hypersurface singularities on $Z_n$. 
Therefore we can obtain a $\mathbb{Q}$-factorialization of it, $W$, by just blowing up $Z_n$ and the exceptional $E_i$. Now consider 
the map $W/T \la X/T$. A $\mathbb{Q}$-factorialization of $Z_0$ can be obtained by running a $(W_0,\epsilon \sum E_i^{W})$ MMP 
over $X_0$. Since $W_0$ has finitely many terminal singularities over any point of $X$, then 
by Propositions 2.11 and 2.12, the operations of the MMP extend in the family. Therefore we obtain 
a morphism $Y/T \stackrel{f}{\la} X/T$, such that $Y_0 \stackrel{f_0}{\la} X_0$ is a $\mathbb{Q}$-factorialization of $Z_0$. 
We now want to conclude that $Y_t \stackrel{f_t}{\la} X_t$ is a $\mathbb{Q}$-factorialization of $Z_t$ for $t$ in a neighborhood 
of $0$. As we have seen in example 2.10, this is not automatic in a family of canonical threefolds. However, the conditions of 
theorem~\ref{deform} mean that $e(X_0)=e(X_t)$, for all $t$. Since $Y_0 \la X_0$ is small, then by semicontinouity of fiber 
dimensions it follows that $Y_t \la X_t$ is also small in a neighborhood of $0$. Therefore, $e(Y_0)=e(Y_t)$. But this means that 
when we run the  $(W_0,\epsilon \sum E_i^{W})$ MMP, it is not possible that there is a divisorial contraction in the central fiber 
and small away from it which must be the case if $Y_t$ fails to be $\mathbb{Q}$-factorial. The theorem now follows.
$\Box$

The next example shows that the conditions of Theorem~\ref{deform} are needed in order for the $D$-minimal models to 
form a family.

\begin{example}\label{ex1}
In this example we construct a one parameter family of canonical 3-folds $Y \la \mathbb{C}^1$ and a family of Weil divisors 
$E$ in $Y$, such that the $E_u$-minimal models of $Y_u$ do not form a family, where $u$ is the parameter.

Let $X_0 \subset \mathbb{C}^4$ be the 3-fold terminal singularity given by $xy+z^3+t^3=0$, where $x,\; y,\; z,\; t$ are the 
coordinates of $\mathbb{C}^4$. Let $X \la \mathbb{C}^1$ be a one parameter deformation of $X_0$ given by 
$uz^2+xy+z^3+t^3=0$, where $u$ is the parameter. $X_0$ has a $cA_2$ singularity but $X_t$ only $cA_1$, for $t\neq 0$.
Let $\Gamma \subset X$ be the plane given by $x=z=t=0$ 
and let $f \colon Y=B_{\Gamma}X\la X$ be the blow up of $X$ along $\Gamma$. Then an explicit calculation of $Y$ shows that 
the exceptional set of $f$ consists of two divisors. A \pone-bundle $E$ over $\Gamma$, and a $\mathbb{P}^2$-bundle $F$ over 
the line $x=y=z=t=0$. Moreover, $E_u$ and $F_u$ are not $\mathbb{Q}$-Cartier for all $u$. 

\noindent\textit{Claim:} The $E_u$-minimal models of $Y_u$ do not form a family. Suppose they do. Then there is a 
morphism $g \colon Z \la Y$ such that $g_u \colon Z_u \la Y_u$ is the $E_u$-minimal model of $Y_u$. Hence $E^Z_u$ and 
$F^Z_u$ are $\mathbb{Q}$-Cartier. Now we can contract $F^Z_0$ to a terminal singularity. Let $p_0 \colon Z_0 \la W_0$ be 
the contraction. Then the whole process is described by the following diagram:
\[
\xymatrix{
    & Z_0 \ar[dl]_{g_0} \ar[dr]^{p_0} &             \\
Y_0 \ar[dr]_{f_0} &  & W_0 \ar[dl]^{q_0} \\
                  & X_0  &}
\]
where $q_0\colon W_0 \la X_0$ contracts $E^W_0$ onto the curve $\Gamma_0 : \; x=z=t=0$. 
If the $E_u$-minimal models of $Y_u$ form a family and a morphism $g \colon Z \la Y$ exists as above, then $p_0$ and $q_0$ 
extend to the family~\cite[Proposition 11.4]{Ko-Mo92}, and the above diagram is a specialization of the diagram 
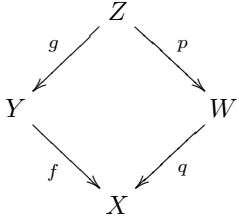
\[
\xymatrix{
    & Z \ar[dl]_{g} \ar[dr]^{p} &             \\
Y \ar[dr]_{f} &  & W \ar[dl]^{q} \\
                  & X  &}
\]
where $W$ is a one parameter deformation of $W_0$.
In particular there is a family of divisorial contractions $q \colon W \la X$ such that for all $u$, $q_u \colon W_u \la X_u$ 
is a divisorial contraction contracting an irreducible divisor $E^W_u$ onto the line $\Gamma_u \subset X_u$ given by $x=z=t=0$. 
However,

\noindent\textit{Claim:} $W_0$ is terminal of index 3, but $W_u$ is terminal of index 2, for all $u \neq 0$.
Since the index is constant in families of terminal 3-folds, 
we get a contradiction and therefore the $E_u$-minimal models of $Y_u$ do not form a family.

We now proceed to justify all the previous claims. By construction, $Y_{u}=B_{\Gamma_{u}}X_{u}$.
In the affine chart given by $x=xt$, $z=zt$, $Y$ is given by \[
u z^2t + xy + z^3t^2+t^2=0.\]
Set $t=0$ to see that $f^{-1}(\Gamma)=E+F$, where $E$ is given by $x=t=0$ and $F$ by $y=t=0$. The other charts are checked similarly. 

Assuming that the $E_u$-minimal models of $Y_u$ form a family, the only thing that needs to be proved is the claim 
about the indices of $W_u$. 
$X_0$ is given by $xy+z^3+t^3=0$ and $\Gamma_0$ by $x=z=t=0$. The general hyperplane section of $X_0$ 
through $\Gamma_0$ has an  $A_2$ DuVal singularity. 
For $u \neq 0$, $X_{u}$ is given by $u z^2 + xy+z^3+t^3=0$. In this case, the general hyperplane section 
through $\Gamma_{u}$ 
has an $A_1$ DuVal singularity. Therefore, $W_{u}$ has index $2$ if $u \neq 0$~\cite[Theorem 5.1]{Tzi02}. 

We will construct the contraction $W_0 \la X_0$ explicitely. $Y_0$ is given by $xy+z^3t^2+t^2=0$. 
$E_0$ is given by $x=t=0$ and $F_0$ by $y=t=0$. 
We want to construct the $E_0$-minimal model of $Y_0$, $Z_0$, and then contract $F_0$ to obtain $W_0$. 
$Y_0$ is given by $xy+(z^3+1)t^2=0$ and a straightforward calculation shows that its singular locus is the line $x=y=t=0$. 
Moreover, $2E_0$ is Cartier at all points except at the three points given by $x=y=t=z^3+1=0$. 
Therefore, in order to make $E_u$ $\mathbb{Q}$-Cartier all we have to do is blow up the ideal $I=(x,t^2)$ which defines $2E_0$. 
Let $g \colon Z_0 \la Y_0$ be the blow up of $I$. Then $g$ is an isomorphism in codimension 1 and moreover 
I claim that the $g$-exceptional curves $C$, are not contained 
in the birational transform $F_0^Z$ and therefore $F_0^Z=F_0=\mathbb{P}^2$. Then $F_0^Z$ can indeed be contracted to a point 
by a $K_{Z_0}$-negative extremal contraction as claimed 
earlier. To see this describe $Z_0$ explicitely. By its construction, $Z_0 \subset \mathbb{C}^4 \times \mathbb{P}^1_{u,v}$ 
and is given by 
\begin{gather*}
xu-t^2v=0 \\
xy+z^3t^2+t^2=0
\end{gather*}
In the affine chart $v=1$, $Z_0$ is given by 
\begin{gather*}
xu-t^2=0 \\
y+z^3u+u=0
\end{gather*}
The $g$-exceptional curves are given by $x=y=t=z^3+1=0$, and $F_0^{Z_0}$ by $y=t=u=0$. Therefore 
no $g$-exceptional curve is contained 
in $F_0^{Z_0}$ and hence $F_0^{Z_0}=F_0=\mathbb{P}^2$. Let $l \subset F_0^Z$ be a line. 
Then $l \cdot K_{Z_0} = l \cdot K_{Y_0}=-1 < 0$ and therefore
$Z_0 \stackrel{p_0}{\la} W_0$ is a $K_{Z_0}$-negative extremal contraction as claimed. 

We now proceed to find the index of $W_0$. We know that $2E_0^{Z_0}$ and $2F_0^{Z_0}$ are Cartier. 
Moreover there is a rational number $a>0$ such that 
\begin{equation}
K_{Z_0}=p^{\ast}K_{W_0}+aF_0^{Z_0}.
\end{equation}
Let $l \subset F_0^{Z_0}=\mathbb{P}^2$ be a line. Then $l \cdot F_0^{Z_0} = -b/2$, where $b$ is a positive integer. 
We want to find $b$. \[
l \cdot (E_0^{Z_0}+F_0^{Z_0})= l \cdot g_0^{\ast}(E_0+F_0)=l \cdot (E_0+F_0)=-1.\]
Now, $l \cdot E_0^{Z_0}=1/2(2E_0^{Z_0}) \cdot l$. This can easily be computed. $ 2E_0^{Z_0}$ is given by $x=t^2=y+z^3u+u=0$ and $F_0^{Z_0}$ 
by $y=t=u=0$. We can take any line in $F_0^{Z_0}$, so take the one given by $y=t=z=u=0$. 
The intersection $l \cap 2E_0^{Z_0}$ is a reduced point and therefore $l \cdot 2E_0^{Z_0} = 1$. Hence \[
l \cdot F_0^{Z_0}=-1-l \cdot E_0^{Z_0}=-3/2.\]
Intersecting $(6)$ with $l$ we see that $a=2/3$. Therefore the index of $W_0$ is 3.
\end{example}

\section{Applications to divisorial contractions.}
In this section we obtain a complete classification of terminal 3-fold divisorial contractions 
$f \colon E \subset Y \la \Gamma \subset X$, 
in the case that $\Gamma$ is a smooth curve and the general hyperplane section $S$ of $X$ through $\Gamma$ has a $D_5$ 
DuVal singularity. 
 
\begin{definition}
A divisorial contraction is a proper morphism $f: E \subset Y \la \Gamma \subset X $, such that $Y$ is $\mathbb{Q}$-factorial,
$Y-E \cong X-\Gamma$, $E$ is an irreducible divisor and  $-K_Y$ is $f$-ample.
\end{definition}
The contraction is called terminal if both $X$ and $Y$ have terminal singularities. We also need  the following definition.
\begin{definition}[\cite{Tzi02}]
Let $0 \in S$ be a $D_n$ DuVal singularity, and let $0 \in \Gamma \subset S$ be a 
smooth curve through the singularity. Let $U \stackrel{f}{\la} S$ 
be the minimal resolution and $\Gamma^{\prime}=f_{\ast}^{-1} \Gamma$. Let \[
\begin{array}{cc}
                                                                 &  \stackrel{E_{n-1}}{\circ} \\
                                                                 & \mid                       \\
\stackrel{E_1}{\circ} \mbox{\noindent ---} \stackrel{E_2}{\circ} \mbox{\noindent ---}  \cdots \mbox{\noindent ---}  & \stackrel{E_{n-2}}{\circ}  \\
                                                                  & \mid \\
                                                                  & \stackrel{E_n}{\circ}
\end{array}. \]
be the dual graph. Then $\Gamma^{\prime}$ intersects either $E_1$, or $E_{n-1}$ or $E_{n}$. In the first case
we call $\Gamma$ of type $DF_l$, and in the others of type $DF_r$.
\end{definition}

\begin{construction} Let $\Gamma \subset X$ be a curve in a 3-fold $X$ 
having at most index 1 terminal singularities. Suppose that the general hyperplane section $S$ of $X$ through $\Gamma$ is 
DuVal and that the curve has at worst plane singularities. Then there is a divisorial contraction $f \colon Y \la X$ contracting 
an irreducible divisor $E$ onto $\Gamma$ such that $Y$ has at worst canonical singularities. Moreover such contraction 
is unique and can be obtained from the following diagram~\cite{Tzi02}.
\[
\xymatrix{ 
Z \ar[d]_h \ar@{-->}[rr]^{\phi} & & Z^{\prime} \ar[d]^{h^{\prime}}\\
W \ar[dr]_g &   & Y \ar[dl]^f \\
            & X & }
\]
$W$ is the blow up of $X$ along $\Gamma$. There are two $g$-exceptional divisors. A ruled surface $E$ over $\Gamma$, and 
a $F=\mathbb{P}^2$ over the singularity. $Z$ is the $E$-minimal model of $W$. After a sequence of flips, $F$ can 
be contracted by a birational morphism $h^{\prime} \colon Z^{\prime} \la Y$, and we obtain the required contraction.  

In order to understand when a terminal contraction exists, i.e., when $Y$ is terminal, it is nececssary to study the singularities 
of $Z$ and $Z^{\prime}$. In general, calculating $D$-minimal models directly is difficult and the possible appearance of flips 
makes things worse. To overcome this difficulty we will use the one parameter theory developed in the first part of the paper 
to degenerate $\Gamma \subset X$ to a $\Gamma_0 \subset X_0$ having simpler equation and therefore easier to manipulate. 
Then we want to deform $Y$ to $Y_0$. In order to do this we must study the following problem.
\end{construction}
\noindent\textbf{Problem:} Let $X \la T$ be a one parameter family of terminal 3-folds and let $\Gamma \subset X$ be flat over $T$, 
such that $\Gamma_t \subset X_t$ is a smooth curve for all $t$. Does there exist a morphism $f \colon Y \la X$ such that 
$f_t \colon Y_t \la X_t$ is a divisorial contraction contracting an irreducible divisor $E_t$ onto $\Gamma_t$?

In general the answer is no. However under certain conditions they do form families, as the next result shows.
\begin{corollary}
Let $X \la T$ be a 1-dimensional family of terminal 3-fold singularities $P_t \in X_t$. Let $\Gamma \subset X$ be flat over $T$
such that $\Gamma_t \subset X_t$ is a smooth curve through $P_t$. Let $H_t$ be the general hyperplane section of $X_t$ 
through $\Gamma_t$. Suppose it is an $A_n$ DuVal singularity for all $t$ and that $\Gamma_t$ intersects an end of 
the fundamental cycle of $H_t$ for all $t$. Then there exist a 
morphism $Y \stackrel{f}{\la} X$, 
such that $Y_t \stackrel{f_t}{\la} X_t$ is a terminal contraction contracting an irreducible divisor $E_t$ onto $\Gamma_t$. 
\end{corollary}
\begin{proof}
Let $g \colon W \la X$ be the blow up of $X$ along $\Gamma$. Then as mentioned earlier, there are two $g$-exceptional divisors 
$E$ and $F$. Moreover, from~\cite{Tzi02} it follows that $W_t$ is $cDV$ for all $t$. 
Moreover, the singular locus of $W_t$ is the line $L_t=E_t\cap F_t$, and $W_t$ is an $A_{n-1}$ DuVal 
singularity at the generic point of $L_t$. Therefore the conditions 
of theorem~\ref{deform} are satisfied and hence there is a morphism $h \colon Z \la W$, such that $Z_0 \la W_0$ is 
the $E_0$-minimal model of $W_0$. Now the construction described above can be done in families and hence there is a 
morphism $f \colon Y \la X$ with the properties claimed. 
\end{proof}
Note that the conditions of Corollary 3.4 are necessary as example 2.16 shows.

We can use the previous result to improve the result of~\cite[Theorem 5.1]{Tzi02}. 
\begin{corollary}
Let $E \subset Y \la \Gamma \subset X$ be a terminal 3-fold divisorial contraction contracting an irreducible divisor $E$ 
onto a smooth curve $\Gamma$. Suppose that the general hyperplane section $S$ of $X$ through $\Gamma$ is $A_n$, and that 
$\Gamma$ intersects an end of the fundamental cycle of $S$. Then $Y$ has index $n+1$.
\end{corollary} 
\begin{proof}
By~\cite{Tzi02} it is possible to write the equation of $X$ as \[
xy+z^{n+1}+t \phi(x,y,z,t)=0,\]
where the curve $\Gamma$ is given by $x=z=t=0$. By a result of Hironaka-Rossi, the equation of $X$ is equivalent to \[
xy+z^{n+1}+t \phi(x,y,z,t) +t^m=0,\]
for $m >>1$. Now we can deform $X$ to $X_0$ given by $xy+z^{n+1}+t^m=0$, and it can be explicitely seen with the same method as in
 example 2.16, that there is a terminal contraction $E_0 \subset Y_0 \la \Gamma_0 \subset X_0$, and $Y_0$ has index $n+1$. 
These contractions form a family by Corollary 3.4.
The index is constant in families of terminal 3-folds and hence the corollary follows.
\end{proof}

\begin{theorem}\label{D5}
Let $0 \in \Gamma \subset S \subset X$. Suppose that $P \in S$ is a $D_{5}$ singularity for the general $S$ through $\Gamma$. Then
\begin{enumerate}
\item If $\Gamma \subset S$ is of type $DF_l$, then there is no terminal contraction.
\item If $\Gamma \subset S$ is of type $DF_r$, then 
\begin{enumerate}
\item If $0 \in X$ is a $cD_4$ singularity, then there is a terminal contraction.
\item If $0 \in X$ is a $cD_5$ singularity, then it is always possible to write the equation of $\Gamma \subset S \subset X$ as \[
x^2+y^2z+xz^2+t[xz\psi(z,t)+axt^k+\phi_{\geq 3}(y,z,t)]=0, \]
so that $z^{\nu}$ does not appear in $\phi_{\geq 3}(y,z,t)$ for any $\nu$, $k \geq 1$ and $I_{\Gamma} =(x,y,t)$. 
Let $a_{i,j,k}$ denote the coefficient of $y^i z^j t^k$ in $\phi_{\geq 3}(y,z,t)$ and $b$ the coefficient of $xt^2$ 
in the above equation (i.e., either 0 or $a$). Then a terminal contraction exists unless 
\begin{enumerate}
\item \[
a_{0,0,4}=a_{1,0,2}=2a_{0,1,2}-b\psi (0,0)=4a_{0,0,3}-b^2=0 \]
or
\item \begin{gather*}
a_{0,2,1}^2-ba_{0,2,1}+a_{0,0,3}=0 \\
a_{0,1,2}-a_{0,2,1}\psi(0,0)=0.
\end{gather*}
\end{enumerate}
\end{enumerate}
\end{enumerate}
If there is a terminal contraction $E\subset Y \stackrel{f}{\la} \Gamma \subset X$, then $Y$ has index 4.
\end{theorem}
\begin{example}
Let $0 \in \Gamma \subset S \subset X$ be given by \[
X: x^2+y^2z+2xz^2+t\phi_{\geq 4}(y,z,t)=0,\]
and $\Gamma$ by $x=y=t=0$. Then there is no terminal divisorial contraction $Y\la X$ contracting an irreducible suface onto $\Gamma$.
\end{example}
\begin{corollary}
Let $\Gamma \subset X$ as in the theorem. Let $E\subset Y \stackrel{f}{\la} \Gamma \subset X$ the canonical divisorial contraction, 
which by~\cite{Tzi02} always exists. Let $\Sigma$ be the general section of $X$ through $0$. 
Then $Y$ is terminal iff $\Sigma^{\prime}=f_{\ast}^{-1}\Sigma$ is normal.
\end{corollary}
The corollary follows immediately from the proof of Theorem~\ref{D5}. At this point I would like to point out that the conclusion 
of the corollary is still true in the case that the general hyperplane section $S$ of $X$ through $\Gamma$ is a $D_{2n}$ DuVal 
singularity as well. However it is not true in general. In particular there are examples when $0 \in X$ is a $cA_2$ singularity 
and $\Sigma^{\prime}$ is not normal~\cite{Tzi03}. 

\begin{proof}[Proof of Theorem~\ref{D5}]
The method that we are going to use is based on the method that appears in~\cite{Tzi02} which was described in the 
beginning of this section. Fix notation as in 3.3.

The idea of the proof is the following. First we obtain normal forms for the equations of $\Gamma$ and $X$. Then we distinguish 
two cases with respect to the equations of $\Gamma \subset X$. A general one and a special one, and we treat them seperately. 
In the general case, we degenerate $\Gamma \subset X$ to a $\Gamma_0 \subset X_0$ whose equation is simpler. Then we follow the method 
described in 3.3 to construct the contraction $Y_0 \la X_0$. We then use the one parameter theory developed in the first part of the 
paper to show that the contractions form a family. Now a deformation of a terminal singularity is again terminal~\cite{Nak98} and 
hence we obtain information about the singularities of $Y$ from what we know about $Y_0$. 

In the special case we work explicitely with the equation of $X$ and $\Gamma$. 
The main difficulty is to describe $Z$ and the possible existence of flips. 
We show that $Z$ can be obtained from the following diagram:
\[
\xymatrix{
    & W_2 \ar[dl]_{h_2} \ar[dr]^{h^{\prime}} & \\
 W_1 \ar[dr]_{h_1} &                                & Z \ar[dl]^h \\
                  & W                             & }
\]
where $W_1$ is the blow up of $W$ along $E$. Let $F_1$ be the $h_1$-exceptional divisor. 
Then by blowing up a suitable multiple of $E^{W_1}$ we make it $\mathbb{Q}$-Cartier and 
get the $E^{W_1}$-minimal model of $W_1$ which we call $W_2$. 
We then show that no $h_2$-exceptional curve is contained in $F_1^{W_2}$ and hence we can 
contract $F_1^{W_2}$ to obtain $Z$. Moreover from the above construction follows that no $h$-exceptional curve is contained 
in $F^Z$, and therefore $F^Z\cong F\cong \mathbb{P}^2$. Hence $F^Z$ can be contracted to a terminal singularity 
and therefore no flips exist in the construction of 3.3. Then in order to decide whether $Y$ is terminal or only canonical 
we need to study the singularities of $Z$ away from $F^Z$. If they are isolated terminal, then so is $Y$. If not, then 
$Y$ is only canonical. 

Part 1. of the theorem follows from~\cite[Theorem 6.1]{Tzi02}. We now proceed to justify all the steps described above. 
By~\cite[Proposition 4.8]{Tzi02}, in suitable analytic coordinates, the equation of $X$ is 
\begin{equation}
x^2+y^2z+xz^2+t[xz\psi(z,t)+axt^k+\phi_{\geq 2}(y,z,t)]=0
\end{equation}
with $k \geq 1$.Moreover, $\Gamma$ is given by $x=y=t=0$, and the monomials $y^2$, $yz$ and $z^{\nu}$ 
do not appear in $\phi_{\geq 2}(y,z,t)$, for any $\nu$.

The appearance of the monomials $yt$ and $t^2$ in $ \phi_{\geq 2}(y,z,t)$ complicate the calculations a lot and it is best 
to consider two cases with respect to these.

\noindent\textit{Case 1.} One or both of the monomials $yt$ and $t^2$ appear in $\phi_2(y,z,t)$. Suppose that $yt$ exists.  
Then write the equation of $X$ as \[
x^2+y^2z+xz^2+yt^2+t\psi(x,y,z,t)=0.\]
We can now deform $X$ to $X_0$ given by \[
x^2+y^2z+xz^2+yt^2=0.\]
Let $\mathbf{X} \la \Delta $ be the deformation over the unit disk. 
Let $\mathbf{\Gamma} \subset \mathbf{X}$ be the deformation of $\Gamma$ given by $x=y=t=0$, and 
let $\mathbf{W} \stackrel{\mathbf{g}}{\la} \mathbf{X}$ be the blow up 
 of $\mathbf{X}$ along $\mathbf{\Gamma}$. Let $\mathbf{E}$, $\mathbf{F}$ be the two $\mathbf{g}$-exceptional divisors. 
Then $\mathbf{W}_u \cong W=B_{\Gamma}X$ for $u \neq 0$ in $\Delta$, and $\mathbf{W}_0\cong W_0=B_{\Gamma_0}X_0$. 
Now one can check that the family $\mathbf{W}$ satisfies the conditions of theorem~\ref{deform}. 
Therefore there exists a morphism $\mathbf{Z} \la \mathbf{W}$, such that $\mathbf{Z}_u \la \mathbf{W}_u$ is the $\mathbf{E}_u$-
minimal model of $\mathbf{W}_u$. Hence the $E$-minimal model of $W$ and the $E_0$-minimal model of $W_0$ form a family. 
Now the process of constructing a divisorial contraction described in 3.3 can be done in families and therefore 
the divisorial contractions $Y\la X$ and $Y_0 \la X_0$ form a family. Therefore there exists a family $\mathbf{Y}\la \Delta$ 
such that $\mathbf{Y}_0\cong Y_0$, and $\mathbf{Y}_u \cong Y$ for all $u \neq 0$ in $\Delta$. 
Hence if $Y_0$ is terminal, then by~\cite{Nak98} 
$Y$ is also terminal. 

At this point I would like to mention 
that all the arguments so far work for the general $D_{2n+1}$ case. However I do not know how to show that $Y$ deforms to 
$Y_0$ for $n \geq 4$. I believe that a more carefull look at theorem~\ref{deform} will treat the general case but the 
amount of calculations involved exceed the value of the result.

We now proceed to show that $Y_0$ is indeed terminal. The main point is to construct $h_0 \colon Z_0 \la W_0$ 
explicitely and to show that 
\begin{enumerate}
\item No $h_0$-exceptional curve is contained in $F^Z$ and therefore $F^Z\cong \mathbb{P}^2$. Then there are no flips involved 
and $F^Z$ is contracted to a terminal singularity, obtaining $Y_0$.
\item $Z_0$ has isolated terminal singularities away from $E_2^Z$ and hence $W_0$ is terminal.
\end{enumerate}
We can construct $Z_0$ as follows. There is a sequence of crepant blow ups \[
W_4 \stackrel{h_4}{\la} W_3 \stackrel{h_3}{\la} W_2 \stackrel{h_2}{\la} W_1 \stackrel{h_1}{\la} W_0 ,\]
where $h_i$ is the blow up of $W_i$ along the birational transform of $E$ in $W_i$. It is a straightforward calculation to check 
that $W_4$ has terminal hypersurface singularities. The point is that at the generic point of $E \cap F$, $W_0$ is an $A_4$ 
DuVal singularity. Then the process described above is just the blow up of a line through an $A_4$ DuVal singularity that intersects 
an edge of the dual graph. This leads to a crepant resolution of the singularity and hence $W_4$ has isolated terminal hypersurface 
singularities. Then in order to make $E^{W_4}$ Cartier, all we need to do is blow it up. So let 
$W_4^{\prime}\stackrel{h_4^{\prime}}{\la} W_4$ be the blow up of $E^{W_4}$. Let $\psi = h_4^{\prime}\circ h_4\circ h_3 \circ h_2 
\circ h_1$.
Then it is not difficult to check that no $h_4^{\prime}$-exceptional curve is contained in $F^{W_4^{\prime}}$. 
Contracting all the $\psi$-exceptional divisors we get the $E$-minimal model of $W_0$, $Z_0$ with all the properties claimed. 
 
\textit{Case 2.} None of the monomials $yt$ and $t^2$ appear in $\phi_{\geq 2}(y,z,t)$. In this case we work explicitely 
with the equation of $X$ and $\Gamma$. We use the normal form for $X$ given in $(7)$.

We start by describing $W=B_{\Gamma}X$.  
In the chart $x=xt$, $y=yt$ it is given by \[
x^2t+y^2tz+xz^2+xtz\psi(z,t)+axt^{k+1}+\phi_{\geq 2}(yt,z,t)=0.\]
Set $t=0$ to see that $E$ is given by $x=t=0$, and $F$ by $z=t=0$. Let $L=E \cap F:x=z=t=0$. 
Now it is easy to see that $W$ is singular along $L$. This is 
what makes this case so different from the $D_{2n}$ cases where $W$ had only one singularity along $L$, and much more difficult to work.
 
As we said before, we try to show that none of the $h_2$-exceptional curves is contained in $F_1^{W_2}$ and also 
study the singularities of $W_2$ away from $F^{W_2} \cup F_1^{W_2}$. 
Moreover, $Z$ can only have isolated $cDV$ point over a $cDV$ point and therefore we restrict our attention 
to what happens over non $cDV$ points. One can also check the other 
charts and see that all the non $cDV$ point are contained in the first one given by $x=xt$, $y=yt$. So it suffices to do 
all our calculations in that chart. 

Let $h_1 :W_1 \la W$ be the blow up of $W$ along $E$. In the chart $x=xt$, $W_1$ is given by 
\begin{equation}
x^2t^2+y^2z+xz^2+xzt\psi(z,t)+axt^{k+1}+\frac{\textstyle 1}{\textstyle t}\phi_{\geq 2}(yt,z,t)=0.
\end{equation}
Let $F_1$ be the $h_1$-exceptional divisor. Setting $t=0$ we see that
\begin{gather}
F_1: z=t=0\\
E^{W_1}: t=y^2+xz+\frac{1}{z}[\frac{1}{t}\phi_{\geq 2}(yt,z,t)]_{t=0}=0.
\end{gather}
Observe that $F^{W_1}$ does not appear in this chart. Now write $\phi_2(y,z,t)=a_3zt$. 
Then $\frac{1}{t}\phi_2(yt,z,t)=a_3z$. Hence $[\frac{1}{t}\phi_{\geq 2}(yt,z,t)]_{t=0}=a_3z$. 
The singular locus of $W_1$ is contained in $F_1 \cap E^{W_1}$ which is the two lines, $l_d$, given by $z=t=y-d=0$, 
with $d$ such that $d^2+a_3=0$. We will now study the singularities of $W_1$ along these lines. The change of variables 
$y \mapsto y+d$ brings the equations of $l_d$ to $y=z=t=0$, and 
\begin{multline*}
Y_1 : x^2t^2+(y+d)^2z+xz^2+xzt\psi(z,t)+axt^{k+1}\\
+\frac{1}{t}\phi_{\geq 2}((y+d)t,z,t)=0.
\end{multline*}
Look what happens along $l_d$ by making the change of variables $x \mapsto x- \delta$. $Y_1$ is given by 
\begin{multline*}
(x- \delta)^2t^2+(y+d)^2z+(x-\delta)z^2+(x-\delta)zt\psi(z,t)\\
+a(x-\delta)t^{k+1}+a_3z+\frac{1}{t}\phi_{\geq 3}((y+d)t,z,t)=0,
\end{multline*}
and this gives 
\begin{multline}
(x-\delta )^2t^2+y^2z+2dyz+(x-\delta)z^2+(x-\delta)zt\psi(z,t)\\
+a(x-\delta)t^{k+1}+\frac{1}{t}\phi_{\geq 3}((y+d)t,z,t)=0.
\end{multline}
It now follows that $W_1$ is singular along these two lines, and smooth away from them. Therefore we want to calculate $W_2$ in 
a neighborhood of these lines. Now it is possible to write \[
\frac{1}{t}\phi_{\geq3}(yt,z,t)=z^2f(y,z)+t\phi(y,z,t) \]
and \[
\phi(y,z,t)=z\phi_1(y,z)+t\phi_2(y,z,t). \]
Then $(8)$ becomes 
\begin{multline}
x^2t^2+y^2z+xz^2+a_3z+xzt\psi(z,t)+axt^{k+1}\\
+z^2f(y,z)+zt\phi_1(y,z)+t^2\phi_2(y,z,t)=0.
\end{multline}
Moreover, 
\begin{gather*}
F_1:z=t=0 \\      
E^{W_1}:t=y^2+xz+a_3+zf(y,z)=0 \\
\cup l_d :y^2+a_3=z=t=0
\end{gather*}     
Now we want to describe $W_2$. It can be obtained by blowing up a suitable multiple of $E^{W_1}$. 
To find it look what happens along $l_d$. The change of variables $y \mapsto y+d$ and $x \mapsto x-\delta$ 
shows that the ideal defining $2E^{W_1}$  \[
I_{2E^{W_1}}=(t^2, y^2+xz+a_3+xt\psi(z,t)+zf(y,z)+t\phi_1(y,z)) \]
is principal at all but finitely many points of $l_d$. Therefore we can get $W_2$ by blowing up $2E^{W_1}$. 
At this point I want to say that the reason of the condition that $yt$ does not exist in the equation of $X$ is 
that if it appears then I don't see which multiple of $E^{W_1}$ has to be blown up to get $W_2$. 
Now $W_2 \subset \mathbb{C}^4_{x,y,z,t} \times \mathbb{P}^{1}_{u,w}$ is given by the equations 
\begin{gather}
wt^2-u[y^2+xz+zf(y,z)+xt\psi(z,t)+t\phi_1(y,z)]=0\\
x^2t^2+y^2z+xz^2+a_3z+xzt\psi(z,t)+axt^{k+1} \nonumber \\
+z^2f(y,z)+zt\phi_1(y,z)+t^2\phi_2(y,z,t)=0.
\end{gather}
In the affine chart $u=1$, $W_2$ is given by 
\begin{gather}
y^2+xz+zf(y,z)+xt\psi(z,t)+t\phi_1(y,z)-wt^2=0\\
\Phi(y,z,t)=x^2+wz+axt^{k-1}+\phi_2(y,z,t)=0.
\end{gather}
Now we study what happens over $l_d$. The curves $C=h_2^{-1}(l_d)$ are given by 
\begin{gather}
y^2+a_3=z=t=0\\
x^2+axt^{k-1}\mid_{t=0} + \phi_2(y,0,0)=0. \nonumber
\end{gather}
From the description of $W_2$ it follows that no component of $C$ is contained in $F_1^{W_2}$. 
Next we want to see what kind of singularities $W_2$ has along $C$. The jacobian, $J$, of 
$W_2$ along $C$ is
\[
\begin{array}{c}
J= \\
\left( \begin{array}{ccccc}
0           &       2d              & x_0 +f(0,0)   & x_0\psi(0,0)+\phi_1(0,0)   & 0     \\
2x_0+at^{k-1}|_{t=0}  & \frac{\textstyle \partial \phi_2}{\textstyle \partial y}(d,0,0) 
&w+\frac{\textstyle \partial \phi_2}{\textstyle \partial z}(d,0,0) & 
\frac{\textstyle \partial \Phi}{\textstyle \partial t}(d,0,0) & 0
\end{array} \right)
\end{array} \]
where $d^2+a_3=0$ and $x_0^2+at^{k-1}|_{t=0}x_0+\phi_2(d,0,0)=0$. If $d \neq 0$ then from $(12)$ it follows that $W_1$ is $cA_n$ 
along $l_d$ and therefore $W_2$ can have at worst isolated cDV points on $C$. Hence the case of interest is when $d=0$ which 
implies that $a_3=0$ as well.
 
Now $x_0\psi(0,0)+\phi_1(0,0)$ is the coefficient of $zt$ in $(12)$ at any point on $l$ given by $x \mapsto x+x_0$. 
Moreover, $x_0^2+at^{k-1}|_{t=0}x_0+\phi_2(0,0,0)$ is the coefficient of $t^2$. 
Hence if it is nonzero then  $W_1$ is $cA_n$ under $C$, and therefore $W_2$ has isolated $cDV$ points along $C$. 
So we want to investigate what happens if $x_0\psi(0,0)+\phi_1(0,0)=0$. 
Now by looking at the Jacobian matrix $J$ we see that $W_2$ is either singular along $C$ or smooth. In the first case 
the resulting contraction will be only canonical and in the second we will have to see what happens in the other chart. 
$W_2$ is singular along $C$ iff 
\begin{gather}
(2x_0+at^{k-1}|_{t=0})(x_0+f(0,0))=0 \\
[(a(k-1)x_0t^{k-2})\mid_{t=0}+\frac{\textstyle \partial \phi_2}{\textstyle \partial t}(0,0,0)]
( x_0+f(0,0))=0 \nonumber \\
\frac{\textstyle \partial \phi_2}{\textstyle \partial y}(0,0,0)(x_0+f(0,0))=0  \nonumber \\
 x_0\psi(0,0)+\phi_1(0,0)=0 \nonumber \\
x_0^2+ax_0t^{k-1}|_{t=0}+\phi_2(0,0,0)=0 \nonumber
\end{gather}

A similar calculation in the other affine piece of $W_2$ given by $w=1$ leads to the same conclusion about the singularities 
of $W_2$ and we again obtain the equations $(18)$. 

Now we will show that there are no other $h_2$-exceptional curves that appear over the other affine piece 
of $W_1$ given by $t=tx$, and therefore the conditions given by $(16)$ are necessary and sufficient. $W_1$ is given by 
\begin{multline}
x^2t+y^2tz+z^2+xzt\psi(z,xt)\\
+ax^{k+1}t^{k+1}+\frac{1}{x}\phi_{\geq 2}(xyt,z,xt)=0.
\end{multline}
Moreover, $x=0$ gives that $F_1:z=x=0$, $F^{W_1}:z=t=0$ and \[
E^{W_1}=(x,y^2t+z+\frac{1}{z}[\frac{1}{x}\phi_{\geq 2}(xyt,z,xt)]_{x=0}). 
\]
The singularities of $W_1$ lie on $L_1 \cup l_d =E^{W_1} \cap F_1$. The line $L_1:x=z=t=0$ is over $L$, and $l_d:x=y-d=t=0$, with 
$d^2+adt^{k-1}|_{t=0}+\phi_2(0,0,0)=0$ as before, is over a point and we already studied the part 
away from zero in the other chart. So we only 
need to see what happens along $L_1$. 

\noindent\textit{Claim:} $E^{W_1}$ is $\mathbb{Q}$-Cartier along $L_1$. 
Therefore there are no $h_2$-exceptional curves over $L_1$ and hence 
the conditions for nonexistence of a terminal contraction are precisely those given by $(18)$. 

We know that $2F^{W_1}$ is Cartier. Consider the index 1 cover of $F^{W_1}$ \[
\pi :\tilde{W_1} \la W_1. \]
We will show that $F_1^{\tilde{W_1}}$ is $\mathbb{Q}$-Cartier over $L_1$ and therefore $F_1$, and hence $E^{W_1}$ too, is 
$\mathbb{Q}$-Cartier along $L_1$. By definition of the index 1 cover~\cite{Ko-Mo98}, \[
\tilde{W_1}=\text{Spec}(\sheaf_{W_2} \oplus I_{F^{W_1},W_1}) \]

It is easy to see now that $\tilde{W_1}$ is given by 
\begin{multline}
x^2+y^2uw+u^2+xuw\psi(uw,xw^2)+bx^{k+1}w^{2k}+\frac{1}{xw^2}\phi_{\geq 2}(xyw^2,uw,xw^2)=0
\end{multline}
and $\pi_1$ is given by $t=w^2$ and $uw=z$. Moreover, \[
F_1^{\tilde{Y_1}}: u=x=0 \]
and therefore, since $t^2$ does not appear in $\phi_{\geq2}(y,z,t)$, it is $\mathbb{Q}$-Cartier along 
$\pi_1^{-1}(L_1): x=u=w=0$. The claim now follows. 

Therefore, a terminal contraction does not exist iff $(18)$ are satisfied. Let $a_{i,j,k}$ be the coefficient of 
$y^i z^j t^k$ in $\phi_{\geq 2}(y,z,t)$. Then $f(0,0)=a_{0,2,1}$, $\frac{\partial \phi_2}{\partial t}(0,0,0)=a_{0,0,4}$, 
$\frac{\partial \phi_2}{\partial y}(0,0,0)=a_{1,0,2}$, $\phi_2(0,0,0)=a_{0,0,1}$ and $\phi_1(0,0)=a_{0,1,2}$. 
Substituting these to $(18)$ we get the conditions claimed by the statement of the theorem.

We now want to give a geometric interpretetion to the existence of the quadratic part $\phi_2(y,z,t)$ in the equation of $X$. 

\noindent\textit{Claim:} $\phi_2(y,z,t)\neq 0$, iff the general hyperplane section of $X$ through $0$ is a $D_4$ DuVal 
singularity. In other words, $0 \in X$ is a $cD_4$ compound DuVal singularity. 

Suppose that $\phi_2(y,z,t)\neq 0$. Then it is possible to deform $X$ to an $X_0$ given by \[
x^2+y^2z+xz^2+t\phi_2(y,z,t)=0.\]
By assumption, the general hyperplane section $S$ of $X$ through $\Gamma$ is a $D_4$ DuVal singularity. Therefore $0\in X$ 
is $cD_n$ for some $n$ (If it was $cA_n$, then the general hyperplane section through $\Gamma$ is $A_m$ for some $m$~\cite{Tzi02}). 
The Milnor number is upper semicontinuous in families and hence if $X_0$ is $cD_4$, then so is $X$. 
Suppose that $t^2 \in \phi_2(y,z,t)$. Then deform $X_0$ to $X_0^{\prime}$ given by $x^2+y^2z+xz^2+t^3=0$. Then the section $z=t$ 
is a $D_4$ DuVal singularity. We may do the same with the other monomial parts of $\phi_2(y,z,t)$ to see that indeed $0\in X$ 
is $cD_4$.

Conversely, suppose that $0\in X$ is $cD_4$. We will use the following property characterising $D_4$ DuVal singularities. 
Let $0 \in T$ be a $D_n$ DuVal singularity. Then by using the Weierstrass preparation theorem, $T$ is given by $x^2+f(z,t)=0$. 
This is $D_4$ iff the cubic part of the previous equation $f_3(z,t)$ is a product of three distinct linear forms. 

So suppose that $\phi_2(y,z,t)=0$ but $0\in X$ is $cD_4$. Eliminating $x$ in $(7)$ we see that the cubic part of the equation of $X$ 
is $y^2z+t\phi_2(y,z,t)$. If $\phi_2(y,z,t)=0$, then this can never be a product of three distinct linear forms.

Now suppose that $Y$ is terminal. We want to find its index. 
Let $S$ be the general hyperplane section of $X$ contaning $\Gamma$ and $S^{\prime}$ its birational transform in $W$. 
Then by~\cite{Tzi02} it follows that $S^{\prime}$ has exacty one singular point which is an $A_{4}$ 
DuVal singularity. So $W$ is an $A_{4}$ DuVal singularity at the generic point of $L$ and hence $5E$ 
and $5F$ are Cartier at the generic point of $L$. 
Arguing as in the proof of~\cite[Theorem 5.1]{Tzi02} we conclude that \[
\text{index}(E^Z)=\text{index}(F^Z)=5. \]
Note that the above arguments show that $Z=Z^{\prime}$ and $F^Z \cong \mathbb{P}^2$. Let $a>0$ such that 
\begin{equation}
K_Z={h^{\prime}}^{\ast}K_Y+aF^Z.
\end{equation}
Let $l \subset F^Z \cong \mathbb{P}^2$ be a general line. 
Then $l \cdot F^Z = - b/5$, with $b$ a positive integer. We want to find $b$. 
\begin{equation}
l \cdot (E^Z+2F^Z)=l \cdot h^{\ast}(E+2F)=-1.
\end{equation}
Moreover, $l \cdot E^Z = 1/5 (5E^Z)\cdot l$. This is easily computable. I claim that the scheme theoretic intersection of 
$5E^Z$ and $F^Z$ is just $3L$. Then for general $l$, $l \cdot E^Z=3/5$. To see the claim it suffices to 
work at the generic point of $L$. Then as we explained earlier, $Z$ is an $A_4$ surface singularity there 
and $E^Z$, $F^Z$ are two lines 
through the singular point. Moreover, by~\cite[Theorem 4.1, Proposition 4.6]{Tzi02}, 
it is possible to write at the generic point of $L$, 
$Z$ as $xy-z^5=0$, $E^Z$ is given by $x-z=y-z^4=0$, and $F^Z$ by $x-z^2=y-z^3=0$. A straightforward calculation shows that 
$\text{length}(5E^Z \cap F^Z)=3$ (at the generic point of $L$), and hence 
$5E^Z \cap F^Z = 3L$. Therefore from $(22)$ it follows that \[
l \cdot F^Z = 1/2(-1-3/5)=-4/5.\]
Now from $(21)$ it follows that $a=5/4$, and therefore $Y$ has index 4.

Finally, the general section $\Sigma$ of $X$ through $0$ is given by $z=ax+by+ct$. 
If one looks carefully through the calculations in the proof of 
theorem~\ref{D5}, it can be seen that the equations $(18)$ are exactly the conditions for $\Sigma$ not to be normal. 
Hence we get corollary 2.4. 
It is unfortunate that the only proof I know of this result is a computational one.
\end{proof}
\begin{corollary}
Let $0 \in \Gamma \subset X$ as above. Let $\Sigma$ be the general hyperplane section of $X$ through $0$. Write $X$ as deformation 
of $\Sigma$, $X \la \Delta$, using $\Gamma$ as the parameter. 
There is a morphism $\Delta \stackrel{\phi}{\la} Def(\Sigma)$ inducing the above deformation. 
Then there is a subset $Z$ of $Def(\Sigma)$, such that, a terminal divisorial contraction exists iff $Im(\phi) \subset Z$.
\end{corollary}
\begin{proof}
Let $\Sigma^{\prime} \stackrel{f}{\la} \Sigma$ be a birational morphism such that $\Sigma^{\prime}$ is normal, 
$-K_{\Sigma^{\prime}}$ is $f$-ample, and $\Sigma^{\prime}$ admits a terminal 1-parameter deformation. Then there is a morphism 
$Def(\Sigma^{\prime}) \stackrel{\phi}{\la} Def(\Sigma)$ \cite[Prop. 11.4]{Ko-Mo92}. Let $Z \subset Def(\Sigma)$ be the union of the 
images of all such morphisms. We may consider $X$ as a deformation of $\Sigma$ along $\Gamma$, say over a base $\Delta$. 
Let $W \stackrel{g}{\la} X$ be the canonical contraction 
that always exists~\cite{Tzi02}. Then $W$ can be considered as a deformation of $\Sigma^{\prime}=g^{\ast}\Sigma$. If $W$ is terminal, 
then we showed that $\Sigma^{\prime}$ is normal. Hence $Im[\Delta \la Def(\Sigma)] \subset Z$. Conversely, if 
$Im[\Delta \la Def(\Sigma)] \subset Z$, then there is a morphism $\Sigma^{\prime} \la \Sigma$ that lifts to a morphism $W \la X$, where 
$W$ is a terminal deformation of $\Sigma^{\prime}$, and is the required contraction.
\end{proof}

\affiliationone{
   Nikolaos Tziolas\\
   Max Planck Institute f\"{u}r Mathematik\\
Vivatsgasse 7, Bonn, 53111, Germany
   \email{tziolas@mpim-bonn.mpg.de}}
\affiliationtwo{~} 
\affiliationthree{%
   Current address:\\
   Department of Mathematics\\
University of Crete\\
Knossos Ave, Heraklion 71409, Greece
   \email{tziolas@math.uoc.gr}}
\end{document}